\documentclass{article}
\usepackage{arxiv}
\usepackage{amsmath,amsfonts,amssymb,amsthm}
\usepackage{graphicx,xcolor}
\usepackage{mathtools}
\usepackage{algorithm}
\usepackage{algorithmic}
\usepackage{scalefnt}
\usepackage{morefloats}
\usepackage{verbatim}
\usepackage{setspace}
\usepackage{multirow}
\usepackage{stmaryrd}
\usepackage{url}  
\usepackage[colorlinks=true]{hyperref}  
\newfont{\notapolice}{cmss8}
\newfont{\rempolice}{cmss9}
\newfont{\tablepolice}{cmtt10}
\newtheoremstyle{note}% 〈name〉
{3pt}% 〈Space above〉
{3pt}% 〈Space below〉
{\notapolice}% 〈Body font〉
{}% 〈Indent amount〉1
{\itshape}% 〈Theorem head font〉
{.}% 〈Punctuation after theorem head〉
{.5em}% 〈Space after theorem head〉2
{}% 〈Theorem head spec (can be left empty, meaning ‘normal’)〉

\extrafloats{100}

\theoremstyle{note}
\newtheorem{rmk}{Remark}
\newtheorem{example}{Example}
\newcommand{\CapText}{\it \small }

\newcommand{\trans}{ {{\textrm t}} }

\newcommand{\balpha}{{\boldsymbol{\alpha}}}
\newcommand{\bbeta}{{\boldsymbol{\beta}}}

\newcommand{\bw}{{\boldsymbol{w}}}

\newcommand{\bc}{{\boldsymbol{c}}}
\newcommand{\bff}{{\boldsymbol{f}}}
\newcommand{\bfg}{{\boldsymbol{g}}}

\newcommand{\bn}{{\boldsymbol{n}}}
\newcommand{\bu}{{\boldsymbol{u}}}
\newcommand{\bx}{{\boldsymbol{x}}}
\newcommand{\bnu}{{\boldsymbol{\nu}}}
\newcommand{\bmu}{{\boldsymbol{\mu}}}
\newcommand{\cloud}{{\mathcal G}}
\newcommand{\incloud}{{\mathring{\mathcal G}}}
\newcommand{\bdcloud}{{\partial \mathcal G}}

\newcommand{\stencil}{{\mathcal{S}}}
\newcommand{\polyset}{{\mathcal{A}}}
\newcommand{\deriveset}{{\mathcal{B}}}

\title{Optimized Stencil Strategy for the Generalized Finite Difference Method: Application to Steady-State Non-Linear Problems}

\author{
     \href{https://orcid.org/0000-0003-2295-5118}{\includegraphics[scale=0.06]{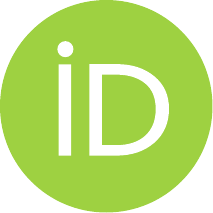}\hspace{1mm}Stéphane Clain} \\
	Centre of Mathematics, Coimbra University, Largo D. Dini, 3000-143 Coimbra, Portugal \\
	\texttt{clain@mat.uc.pt} \\
	\And
	\href{https://orcid.org/0000-0001-9036-1309}{\includegraphics[scale=0.06]{orcid.pdf}\hspace{1mm}Jorge Figueiredo} \\
	Center of Mathematics, Campus de Gualtar, 4710-057 Braga \\
	\texttt{jmfiguei@math.uminho.pt} \\
}

\renewcommand{\headeright}{Article}
\renewcommand{\undertitle}{Article}
\renewcommand{\shorttitle}{Structural schemes for Hamiltonian system}

\hypersetup{
pdftitle={Optimized Stencil Strategy for the Generalized Finite Difference Methodm},
pdfsubject={math.NA},
pdfauthor={Stephane Clain, Jorge Figueiredo},
pdfkeywords={Very high order, Generalised Finite Difference Method, nonlinear PDE},
}

\begin{document}
\maketitle
\begin{abstract}
We propose an optimized stencil strategy for the Generalized Finite Difference Method (GFDM) applied to non-linear problems.
We take advantage of the flexibility of GFDM to engineer specific stencils by agglomerating nodes and balancing size with numerical accuracy. Previous work, focusing on the stencil construction for the linear convection-diffusion problem, showed that optimizing stencils and scaling parameters improves the scheme. The present study aims to investigate similar benefits for non-linear problems such as the Burgers' equations and the weakly compressible Navier-Stokes system.
\end{abstract}

\section{Introduction}
Meshless methods, {\it i.e.} numerical methods based on a set of points with no specific connectivity, have received great attention over the last two decades due to their ability to handle complex geometries, ability to adapt the points density in function of the local regularity and to support both Lagrangian and Eulerian formulations. Such an approach has been broken down into several methods such as Smooth Particles Hydrodynamics (SPH) \cite{GM77,L77}, Finite Pointset method (FPsM) \cite{TK02}, Finite Point method (FPtM) \cite{OI96,OI98,OO13}, or Generalized Finite Difference Method (GFDM) \cite{BUG01,ZL22,WL25}. 

Finite Pointset and Finite Point methods use the entire neighborhood of points belonging to a compact ball, controlled by a threshold parameter included in the weighted function. Moreover, both rely on a moving least-squares method where the reconstruction is achieved at any point (not necessarily a point of the cloud). On the contrary, the Generalized Finite Difference Method depends on particular stencils (not all the points of the neighborhood are necessarily used), and a classical least-squares technique centered on the cloud points is used to derive relations between derivatives and the function values over the stencil. Consequently, the GFDM provides an additional advantage over FPsM and FPtM, the ability to play with the stencil choice and, by that way, obtain the best trade-off between stencil size and accuracy.

In \cite{CF24}, the authors take advantage of several properties inherent to the local representation in the GFDM. The choice of the stencil together with the parameterization of the weight function and the polynomial's scaling parameter enable to optimize the quality of the representation of the derivatives while preserving small-sized stencils. The study was developed for the linear convection diffusion operator in two-dimensional space, where computational cost reduction and better accuracy for a given cloud of points are reported. The question then arises as to whether such gains are still present for non-linear problems. To this end, we conduct a new study on two very popular systems, the viscous Burgers' equation \cite{FL14} and the weakly compressible Navier-Stokes equations (Tait equation of state). We do not integrate any upwind techniques \cite{LF17,R23} since we carry out the simulations with low Reynolds numbers.

\section{Design of the numerical method}
Let $\Omega$ be a bounded open set of $\mathbb R^2$ with boundary $\partial \Omega$ and denote by $\bx=(x_1,x_2)$ a generic point of $\overline \Omega$. We denote by $\phi$ any smooth function defined on $\overline \Omega$ and sufficiently regular such that the first- and second-order derivatives are continuously bounded on $\overline \Omega$. We restrict the presentation to the two-dimensional geometries for the sake of simplicity, but extension to higher dimensions is straightforward.

\subsection{Cloud of points, approximations, and stencil}
A cloud consists of a finite set of points $\bx_i\in\overline \Omega$ characterized by the indices set $\cloud=\{1,\dots,I\}$ which we split into 
$\incloud$, the indices of the points that belong to the open set $\Omega$ and $\bdcloud$ the indices of points on the boundary $\partial \Omega$. For any node $i\in\cloud$, we denote by $\stencil_i$ the associated stencil constituted by a sublist of indices from $\cloud$ such that $i\in\stencil_i$, while $|\stencil_i|$ stands for the stencil size.

We introduce the multi-index notation $\balpha=(\alpha_1,\alpha_2)$, with $\alpha_1,\alpha_2\in \mathbb N_0$ with $|\balpha|=\alpha_1+\alpha_2$ and $\balpha!=(\alpha_1!)(\alpha_2!)$, while $\bx^\balpha=(x_1)^{\alpha_1}(x_2)^{\alpha_2}$. For $i\in\cloud$ we also define the monomial function of degree $\balpha$ centred at point $\bx_i$ with
$$
\pi_\balpha(\bx;\bx_i,h_i)=\left ( \frac{\bx-\bx_i}{h_i}\right )^\balpha
$$
where $h_i>0$ is a scaling parameter.

Given $\bbeta=(\beta_1,\beta_2)$, with $\beta_1,\beta_2\in \mathbb N_0$, we associate the partial derivative $\partial^\bbeta \phi=\partial^{\beta_1}_{x_1}\,\partial^{\beta_2}_{x_2}\phi(x_1,x_2)$. Inequality $\balpha\geq \bbeta$ means that $\alpha_1\geq \beta_1$ and $\alpha_2\geq \beta_2$.
In particular, if $\balpha\geq \bbeta$, one has  
$$
\partial^\bbeta \pi_\balpha(\bx;\bx_i,h_i)
%=\frac{\balpha!}{(\balpha-\bbeta)!}\,\frac{1}{(h_i)^{|\bbeta|}}\, \left ( \frac{\bx-\bx_i}{h_i}\right )^{\balpha-\bbeta}
=\frac{\balpha!}{(\balpha-\bbeta)!}\,\frac{1}{(h_i)^{|\bbeta|}}\,\pi_{\balpha-\bbeta}(\bx;\bx_i,h_i),
$$
while $\partial^\bbeta \pi_\balpha(\bx;\bx_i,h_i)=0$, otherwise. Note that taking $\bx=\bx_i$ provides the relations
\begin{equation}\label{eq::alpha_derivative_monomial}
\partial^\balpha \pi_\balpha(\bx_i;\bx_i,h_i)=\frac{\balpha!}{(h_i)^{|\balpha|}}  \quad\textrm{ and }\quad
\partial^\bbeta \pi_\alpha(\bx_i;\bx_i,h_i)=0, \textrm{ otherwise}.
\end{equation}

\subsection{Polynomial degree set and derivative order set}
The polynomial degree set $\mathcal A$. The space $\mathbb P_i^\polyset$ is constituted by the span of the monomial family
$$
\mathbb P_i^\polyset=\textrm{span}\big\{\pi_\balpha(\bx;x_i,h_i),\ \balpha \in \polyset\big \}
$$
we shall call the $\mathcal A$-polynomial space centred at node $i$.
\begin{rmk}
For the particular case where $\polyset=\{\balpha,\ |\balpha|\leq p\}$, the set simply reads $\mathbb P^p$ since the reference point has no longer importance. Indeed, the centred basis is equivalent to the canonical one. For example, for the polynomial basis of degree $p=2$, we have a basis attached to node $i$
$$
\mathbb P^2=\textrm{span}\Big \{\pi_\alpha(\bx;\bx_i,h_i),\  \balpha\in \polyset=\big \{(0,0),(1,0),(0,1),(2,0),(1,1),(0,2)\big \} \Big\}. 
$$
In that case, the reference node has no importance  since we can rewrite the polynomial space with the canonical monomial. $\blacksquare$
\end{rmk}
Consequently, any polynomial function of $\mathbb P_i^{\polyset}$ is written under the form
$$
\pi(\bx,\bc_i;\bx_i)=\pi(\bx,\bc_i;\bx_i,h_i,\polyset)=\sum_{\balpha\in\polyset} c_{i,\balpha}\pi_\balpha(\bx;\bx_i,h_i),
$$
where the vector $\bc_i$ collects the entries $\big (c_{i,\balpha}\big)_{\balpha\in\polyset}$.
\vskip 1em
The set $\deriveset$ contains the list of derivatives $\bbeta$ we aim at discretizing via  differential operators. For example,
$\big \{\partial^0_{x_1}\partial^0_{x_2},\, \partial^1_{x_1}\partial^0_{x_2},\, \partial^0_{x_1}\partial^1_{x_2},\,\partial^1_{x_1}\partial^1_{x_2}\big \}\iff \deriveset=\big \{(0,0),(1,0),(0,1),(1,1)\big \}$.

Finally, for any regular function $\phi$, we denote  $\phi_i\approx  \phi(\bx_i)$ an approximation of $\phi$ at point $\bx_i$ and for any partial derivative $\partial^\bbeta$, we use the notation  $\phi_i^\bbeta\approx \partial^\bbeta \phi(\bx_i)$. We then organize the data following two vectors:
\begin{itemize}
\item $\Phi_{\stencil_i}=(\phi_j)_{j\in\stencil_i}$~the approximation vector associated with the stencil $\stencil_i$;
\item $\Phi=\Phi_{\cloud}=(\phi_i)_{i\in\cloud}$ the vector that gathers all the approximations associated to the cloud $\cloud$.
\end{itemize}

\subsection{Interpolation and differential operators' discretization}

Given a stencil $\stencil_i$ associated to a node $i\in\cloud$, $\Phi_{\stencil_i}$ the vector of approximations over the stencil, and $\mathbb P^\polyset_i$ the polynomial subspace, we define the following functional for any $\pi\in\mathbb P_i^{\mathcal A}$,
$$
E(\bc_i;\stencil_i,\Phi_{\stencil_i},\bw_i,\polyset,h_i)=\sum_{j\in\stencil_i} w_{ij}\Big [\pi(\bx_j,\bc_i;\bx_i,h_i,\polyset)-\phi_j\Big ]^2,
$$
where  $\bw_i=(w_{ij})_{j\in\stencil_i}$ contains the weights associated to each point of the stencil $\stencil_i$. A usual choice is $w_{ij}=\rho(\Vert \bx_i-\bx_j\Vert;\tau_i)$, where $\rho(s;\tau)$ is a decreasing function of $s\in [0,+\infty[$ controlled by a scaling parameter $\tau$. For instance, the Gaussian weights read
$$
w_{ij}=\rho(\Vert \bx_j-\bx_i\Vert;\tau_i)=\exp \left( - \frac{\Vert \bx_j-\bx_i\Vert^2 }{\tau^2_i} \right).
$$

We define the $|\polyset|\times|\stencil_i|$ matrix $\displaystyle M_i=\big (m_{i;\balpha,j}\big )^{\balpha\in\polyset}_{j\in\stencil_i}$ with
$$
m_{i;\balpha,j}=\pi_\balpha(\bx_j;\bx_i,h_i)=  \left ( \frac{\bx_j-\bx_i}{h_i}\right )^\balpha,\quad \balpha\in\polyset,\, j\in\stencil_i.
$$
Using the Euclidean norm, the functional rewrites under the form
$$
E_i(\bc_i;\stencil_i,\Phi_{\stencil_i},\bw_i,\polyset,h_i)= \Vert W_i(M_i\bc_i-\Phi_{\stencil_i})\Vert^2,
$$
with $W_i=\textrm{diag}(\bw_{i})$ the $|\stencil_i|\times |\stencil_i|$ diagonal matrix of the weights.

Minimization of the functional with respect to vector $\bc_i$ leads to the linear system
$$
(M_i)^\trans W_i M_i\, \bc_i=(M_i)^\trans W_i\Phi_{\stencil_i}.
$$
Assuming that the choice of $\polyset$ and $\stencil_i$ provides a matrix $M_i$ that satisfies the maximal rank property, then $M_i^{\dag}=(M_i)^\trans W_i M_i$ is a non-singular $|\polyset|\times|\polyset|$ matrix and the polynomial coefficients are given by $\bc_i=M_i^\ddag \Phi_{\stencil_i}$ with
$$
M_i^\ddag =\big (M_i^\dag\big  )^{-1}(M_i)^\trans W_i.
$$
Consequently, we deduce the coefficient of a $\mathcal A$-polynomial approximation given by
$$
\widetilde \phi_i(\bx)=\phi_i(\bx,\bc_i;\bx_i,\stencil_i,h_i,\mathcal A)=\sum_{\alpha\in \mathcal A}c_{i,\balpha} \pi_\balpha(\bx;\bx_i,h_i).
$$
Under the condition $\mathcal B\subset \mathcal A$, we get
$$
\partial^\bbeta\pi_\alpha(\bx_i)=\delta_{\balpha,\bbeta}\frac{\balpha!}{(h_i)^\balpha}.
$$
and deduce
\begin{equation}\label{approximation_discretization}
\partial^\beta \widetilde \phi_i(\bx_i)=\frac{\bbeta!}{(h_i)^\bbeta}c_{i,\bbeta}
=\frac{\bbeta !}{(h_i)^{|\bbeta|}}\, \sum_{j\in\stencil_i} m^\ddag_{i;\bbeta,j} \phi_j
=\sum_{j\in\stencil_i} a_{i;\bbeta,j} \phi_j, \textrm{ with } 
a_{i;\bbeta,j}=\frac{\bbeta !}{(h_i)^{|\bbeta|}}\,m^\ddag_{i;\bbeta,j}.
\end{equation}
Relation \eqref{approximation_discretization} provides the discretization of the partial derivatives, in terms of the stencil-node values, which we use in the physical equations.
\begin{rmk}[Pointset method versus generalized finite difference method]
There exists a subtle difference between the two methods when dealing with the reconstruction of $\phi_i(x_i)$. Indeed, we have two ways.
We take $\phi_i(x_i)=\phi_i$ (Generalized Finite Difference Method) or we take $\phi_i(x_i)=\sum_{j\in\stencil_i} a_{i;(0,0),j} \phi_j$ (Pointset Diffuse approximation). The rest of the derivatives' approximation for $\bbeta\neq(0,0)$ are all the same.
\end{rmk}
\begin{example}
Consider the equation $-\Delta \phi=f$. Then the discretization at node $i$ simply reads
$$
-\sum_{j\in\stencil_i} (a_{i;(2,0),j}+a_{i;(2,0),j})\phi_j=f_i.
$$
Similarly, the boundary condition $\nabla \phi\cdot \bn=g$ at node $i$ reads
$$
\sum_{j\in\stencil_i} (a_{i;(1,0),j}n_{i,1}+a_{i;(0,1),j}n_{i,2})\phi_j=g_i.
$$
\end{example}

\subsection{Discretization optimization}

Given a point $\bx_i\in\cloud$, matrix $M_i^\ddag=M_i^\ddag(h_i,\tau_i,\stencil_i)$ depends on the polynomial scaling factor $h_i$, the weights scaling factor $\tau_i$, and the stencil $\stencil_i$. In \cite{CF24} we develop an optimization procedure based on the condition number of matrix $M_i^\dag$ since it corresponds to the non-singular matrix we have to inverse to determine the derivatives. We then introduce a quantification of the polynomial reconstruction efficiency by computing the condition number $\chi_2(M^\dag_i)$ and we aim at determining the parameters that minimize it. In other words, we seek 
$$
(h_i^\star,\tau_i^\star,\stencil^\star_i)=\arg \min \, \chi_2(M_i^\dag(h_i,\tau_i,\stencil_i))
$$
under the constraint $|\stencil_i|\geq |\balpha|$, since we need more points ({\it i.e.}, equations) than unknown polynomial coefficients.

Optimization over a discrete set $\stencil$ and continuous sets $]0,+\infty[\times ]0,+\infty[$ is a difficult task, so the authors combine a heuristic approach for the stencil choice along with a gradient method for the continuous parameters. In this way, we obtain a rather small stencil ($|\polyset|\leq |\stencil_i|\leq 1.5 |\polyset|$) with a dramatic cut in the condition number in comparison with an unsophisticated stencil consisting of the closest points to the reference point $\bx_i$ (about twice the size of $\polyset$ to produce a non-singular matrix $M_i^\dagger$).  
\begin{rmk}
In our previous work the stencil associated with each reference point used to start the heuristic approach is constructed by selecting the 2$|\polyset|$ closest neighbour points. For linear PDEs the strategy enables us to achieve optimal convergence rates, but for certain non-linear problems it occasionally results in suboptimal convergence rates. To improve robustness, we select the 3$|\polyset|$ closest points and discard $|\polyset|$ of them to avoid poor angular distribution (e.g., near collinearity with respect to the reference point). The resulting set of 2$|\polyset|$ points has a more uniform angular spread, and this enriched initial configuration is then used as the input for the subsequent heuristic stencil-optimization procedure.
\end{rmk}
Since we report excellent orders of accuracy for linear PDE problems on complex domains (curved boundary domains) in  \cite{CF24}, we aim to extend the same framework to non-linear PDE problems and systems by benchmarking several important non-linear problems and assessing its efficiency in terms of both accuracy and stability.

\section{The two-dimensional Burgers' equations}
The main goal is to assess the accuracy improvement achieved by the optimized stencil, as well as the robustness of the method. The linear case has been studied in \cite{CF24} where we showed the advantages and gains of the optimization, but the question arises when we are dealing with non-linear problems.

\subsection{Modeling and linearization}
We consider the two-dimensional stationary forced Burgers' equations with viscous term, 
\begin{align*}
u_1\,\partial_{x_1}u_1+u_2\,\partial_{x_2}u_1-\varepsilon \big( \partial_{x_1}^{2}u_1+\partial_{x_2}^{2}u_1\big) & =f_{1},\quad \textrm{ in } \Omega , \\
u_1\,\partial_{x_1}u_2+u_2\,\partial_{x_2}u_2-\varepsilon \big( \partial_{x_1}^{2}u_2+\partial_{x_2}^{2}u_2\big) & =f_{2},\quad \textrm{ in } \Omega ,
\end{align*}
where $\bu=(u_1,u_2)$ is the unknown velocity, $\bff=(f_1,f_2)$ the source term, and $\varepsilon > 0$ the viscosity parameter.
The general Robin boundary conditions are given by
$$
\mu_1(\bx_i)\, u_1+\nu_1(\bx_i) \nabla u_1\cdot \bn=g_1(\bx_i),\quad
\mu_2(\bx_i)\, u_2+\nu_2(\bx_i) \nabla u_2\cdot \bn=g_2(\bx_i),
\quad \textrm{ on } \partial \Omega,
$$
where $\bmu=(\mu_1,\mu_2)$, $\bnu=(\nu_1,\nu_2)$, and $\bfg=(g_1,g_2)$ are given functions and $\bn=\bn(\bx)$ is the outward unit normal vector. Notice that the velocity components are not coupled on the boundary. In this context, we define for any node $i\in\cloud$ the label $\ell(i)$ that indicates the physical equation associated with the node. For example, $\ell(i)=1$ means that one has to satisfy the inner equations, while $\ell(i)=2$ corresponds to the case where the boundary equations hold.

To solve the non-linear system, we use a simple fixed-point method based on a local linearization where we build a sequence $\bu^{[k]}=(u_1^{[k]},u_2^{[k]})$, $k=0,1,\dots$, of approximations of the solution $\bu$ on $\overline \Omega$. 
To this end, assuming $(u_1^{[k]},u_2^{[k]})$ is known, we seek a solution $(u_1^{[k+1]},u_2^{[k+1]})$ of the linearized problem  
\begin{align*}
u_1^{[k]}\,\partial_{x_1}u_1^{[k+1]}+u_2^{[k]}\,\partial_{x_2}u_1^{[k+1]}-\varepsilon\Big (\partial_{x_1}^{2}u_1^{[k+1]}+\partial_{x_2}^{2}u_1^{[k+1]}\Big )&=f_{1}, \\
u_1^{[k]}\,\partial_{x_1}u_2^{[k+1]}+u_2^{[k]}\,\partial_{x_2}u_2^{[k+1]}-\varepsilon\Big (\partial_{x_1}^{2}u_2^{[k+1]}+\partial_{x_2}^{2}u_2^{[k+1]}\Big )&=f_{2}.
\end{align*}
Note that the problem is decoupled and one solves for the velocities $u_1$ and $u_2$ independently. The iterative procedure stops when both $\Vert u_1^{[k+1]}-u_1^{[k]}\Vert$ and $\Vert u_2^{[k+1]}-u_2^{[k]}\Vert$ are lower than a given tolerance. Convergence of the fixed-point is not guaranteed, in particular for large Péclet numbers, since the stencil choice does not take into account the upwind direction, leading to a centered scheme for the convective term. Therefore, benchmarks have been carried out with a relatively low Péclet number to focus on the accuracy and convergence order of the numerical scheme. 

\subsection{Discretization and assembly}
We denote by $U_1$ and $U_2$ the vectors of $\mathbb R^I$ with the numerical solution of the non-linear system, that is, 
$$
U_1=(u_{1,i})_{i\in\cloud}, \quad U_2=(u_{2,i})_{i\in\cloud}.
$$
Moreover, for any point $\bx_i$, $i\in\cloud$, and $\phi=u_1,u_2$, we have $\phi(\bx_i)\approx\phi_i^{(0,0)}=\phi_i$, and we use relations \eqref{approximation_discretization} for the discrete derivatives, that is
\begin{eqnarray}
&\displaystyle \partial_{x_1}\phi(\bx_i)\approx\phi_i^{(1,0)}= \sum_{j\in\stencil_i} a_{i;(1,0),j} \phi_j,
&\displaystyle \partial_{x_2}\phi(\bx_i)\approx\phi_i^{(0,1)}= \sum_{j\in\stencil_i} a_{i;(0,1),j} \phi_j,\label{eq:first_derivative}\\
&\displaystyle \partial^2_{x_1}\phi(\bx_i)\approx\phi_i^{(2,0)}= \sum_{j\in\stencil_i} a_{i;(2,0),j} \phi_j,
&\displaystyle \partial^2_{x_2}\phi(\bx_i)\approx\phi_i^{(0,2)}= \sum_{j\in\stencil_i} a_{i;(0,2),j} \phi_j.\label{eq:second_derivative}
\end{eqnarray}
Since the simulation is carried out using a simple fixed-point method, we construct a sequence of approximations using the following iterative procedure: given ($U_1^{[k]},U_2^{[k]}$), we compute 
$(U_1^{[k+1]},U_2^{[k+1]})$, the approximation of the associated continuous linearized problem at iteration $k+1$:
\begin{align*}
u_{1,i}^{[k],(0,0)}\,u_{1,i}^{[k+1],(1,0)}+u_{2,i}^{[k],(0,0)}\,u_{1,i}^{[k+1],(0,1)}%
-\varepsilon \left( u_{1,i}^{[k+1],(2,0)}+u_{1,i}^{[k+1],(0,2)}\right ) & =f_{1}(\bx_i), \\
u_{1,i}^{[k],(0,0)}\,u_{2,i}^{[k+1],(1,0)}+u_{2,i}^{[k],(0,0)}\,u_{2,i}^{[k+1],(0,1)}%
-\varepsilon \left(u_{2,i}^{[k+1],(2,0)}+u_{2,i}^{[k+1],(0,2)}\right ) & =f_{2}(\bx_i). 
\end{align*}
%\textcolor{red}{where $\phi_i^{[k],\bbeta} \approx \partial^\bbeta\phi(\bx_i)$ at the stage $k$}.
Replacing the derivatives by the local representation given by \eqref{eq:first_derivative} and \eqref{eq:second_derivative},
we then assemble two independent linear systems $A^{[k]}_1U_1=b_1$ and $A^{[k]}_2U_2=b_2$. For example, the assembled matrix corresponding to vector $U_1$ is obtained as follows: 
\begin{itemize}
\item if $\ell(i)=1$, then for $j\in\stencil_i$ we set
$$
A^{[k]}_1[i,j]=u_{1,j}^{[k]}\,a_{i;(1,0),j}+u_{2,j}^{[k]}\,a_{i;(0,1),j}-\varepsilon\Big(a_{i;(2,0),j}+a_{i;(0,2),j} \Big )
$$
and $A^{[k]}_1[i,j]=0$ otherwise. The right-hand side is simply $b_1[i]=f_1(\bx_i)$.
\item if $\ell(i)=2$, then for $j\in\stencil_i$ the discrete Robin condition provides the following entries 
$$
A^{[k]}_1[i,j]=\delta_{i,j}\mu_1(\bx_i)+\nu_1(\bx_i) (a_{i;(1,0),j}n_{i;1}+a_{i;(0,1),j}n_{i;2}),
$$
and $A^{[k]}_1[i,j]=0$ otherwise. The right-hand side is given by $b_1[i]=g_1(\bx_i)$.
Note that whenever the Dirichlet conditions hold, one has $A^{[k]}_1[i,j]=\delta_{i,j}$ and $b_1[i]=u_{D_1}(\bx_i)$.
\end{itemize}
We update the numerical approximation by computing $U_1^{[k+1]}$ and $U_2^{[k+1]}$, solutions of $A^{[k]}_1\,U_1^{[k+1]}=b_1$ and $A^{[k]}_2\,U_2^{[k+1]}=b_2$, respectively.  

%$Re=1/\varepsilon$ is the Reynolds number and 
\subsection{Test case 3.1}
\begin{figure}[ht]
\centering
\includegraphics[width=0.24\textwidth,clip=true,viewport=150 280 485 558]{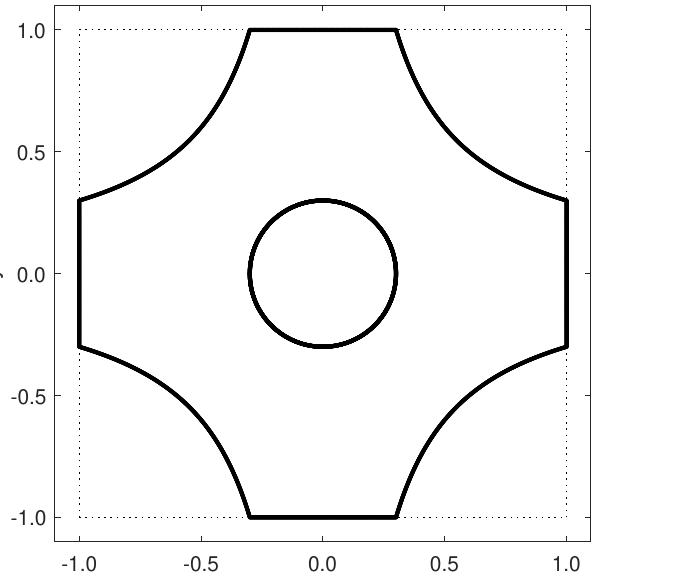}
\includegraphics[width=0.24\textwidth,clip=true,viewport=135 280 470 558]{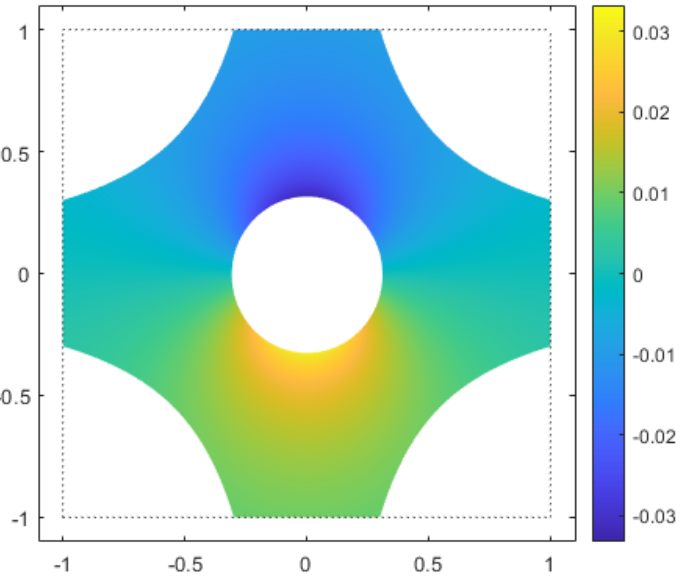}
\includegraphics[width=0.24\textwidth,clip=true,viewport=135 280 470 558]{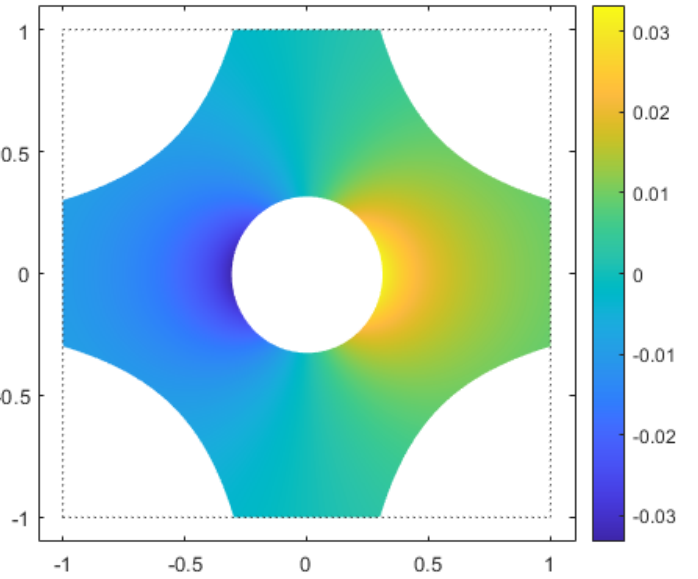}
\includegraphics[width=0.24\textwidth,clip=true,viewport=150 280 485 558]{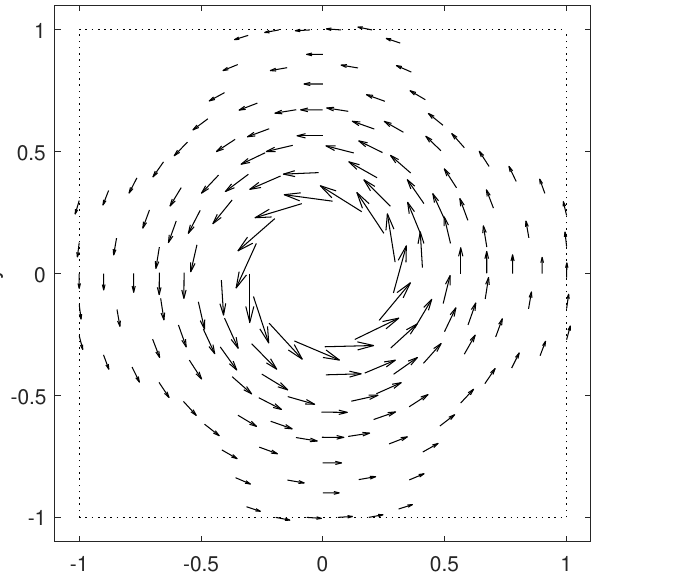}
\caption{\CapText Test case 3.1. Geometrical setting (left panel) and exact solution (from left to right: $u_1$, $u_2$, $\bu$).}
\label{fig:test100solution}
\end{figure}

In this first test case devoted to the two-dimensional Burgers' equations, we prescribe the Dirichlet boundary conditions for both horizontal and vertical velocities, that is $\bmu(\bx_i)=\mathbf{1}$, $\bnu(\bx_i)=\mathbf{0}$, for all $i\in\bdcloud$. We further take $\varepsilon=10^{-2}$ and consider the following irregular domain,
$$
\Omega =\left\{ \bx\in \left] -1,1\right[ \times \left] -1,1\right[
:x_1^{2}+x_2^{2}>0.3^{2}\wedge |x_1x_2|<0.3\right\}.
$$
In order to perform a convergence analysis we set $\bff=(f_1,f_2)$ and $\bu_D=(u_{D_1},u_{D_2})$ such that 
$$
u_1(\bx)=-\frac{\varepsilon x_2}{x_1^{2}+x_2^{2}},\qquad u_2(\bx)=\frac{\varepsilon x_1}{%
x_1^{2}+x_2^{2}}
$$
is the unique solution of the boundary value problem (BVP). We note that this solution is invariant to rotation with respect to the origin of the coordinate
system. The resulting expression for $(f_1,f_2)$ is
$$
f_{1}(\bx)=-\frac{\varepsilon ^{2}x_1}{\left( x_1^{2}+x_2^{2}\right) ^{2}},\qquad
f_{2}(\bx)= \frac{\varepsilon ^{2}x_2}{\left( x_1^{2}+x_2^{2}\right) ^{2}}.
$$
The initial velocity for the fixed-point method, ($U_1^{[0]},U_2^{[0]}$), is taken equal to the average value over the nodes located on $\partial \Omega$. The geometrical setting and the exact solution of the BVP are illustrated in Figure \ref{fig:test100solution}.

\begin{figure}[ht]
\centering
\includegraphics[width=0.24\textwidth,clip=true,viewport=155 280 465 580]{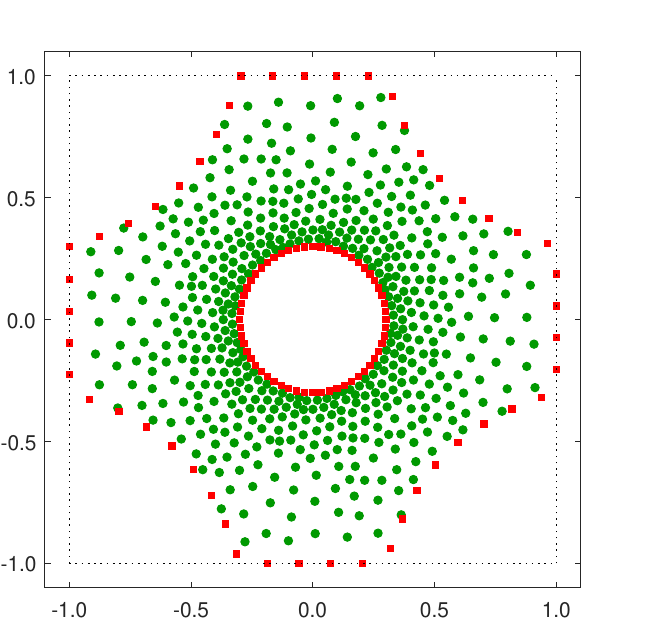}
\includegraphics[width=0.24\textwidth,clip=true,viewport=155 280 465 580]{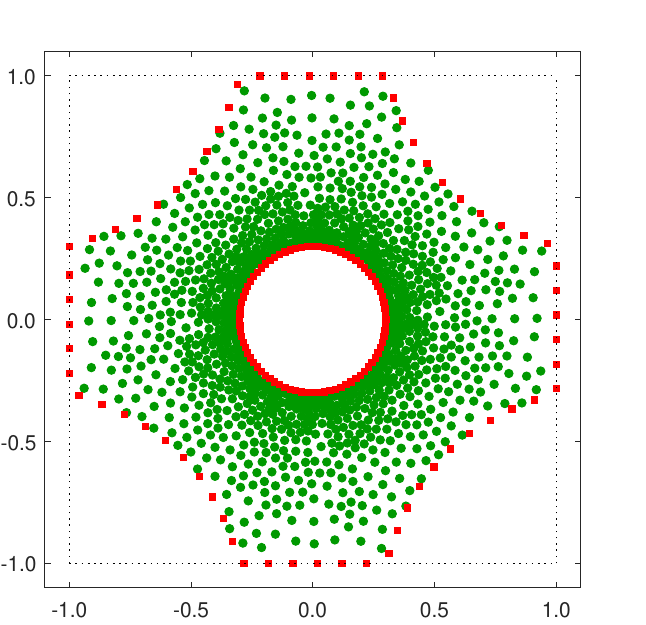}
\includegraphics[width=0.24\textwidth,clip=true,viewport=155 280 465 580]{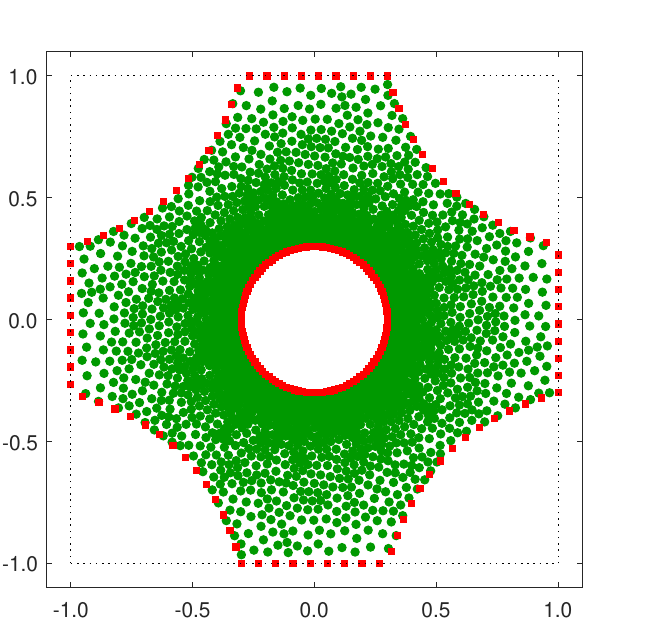}
\includegraphics[width=0.24\textwidth,clip=true,viewport=155 280 465 580]{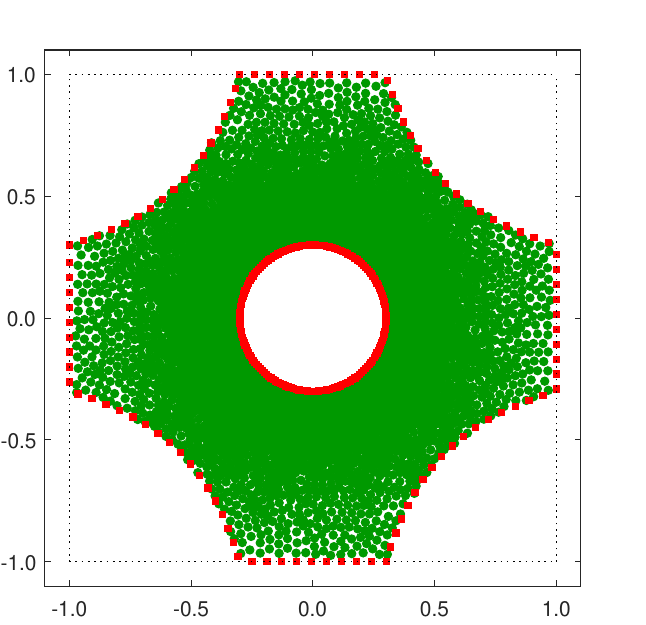}
\caption{\CapText Test case 3.1. For $\overline{\Omega}$ the panels illustrate, from left to right, four clouds obtained with a Delaunay triangulation procedure with $622$, $1288$, $2458$, and $4886$ nodes. The green filled circles and the red squares correspond, respectively, to interior nodes and boundary nodes where Dirichlet conditions are imposed.}
\label{fig:testcase100grid}
\end{figure}

Numerical simulations are carried out with $622$, $1228$, $2458$, and $4886$ nodes, obtained with a Delaunay triangulation procedure using the Gmsh$\textsuperscript{\copyright}$ package. We note that the number of nodes on the inner and outer boundaries is almost the same. Figure \ref{fig:testcase100grid} illustrates the different clouds used in the present test case.

\begin{figure}[ht]
\centering
\includegraphics[width=0.24\textwidth,clip=true,viewport=165 265 465 550]{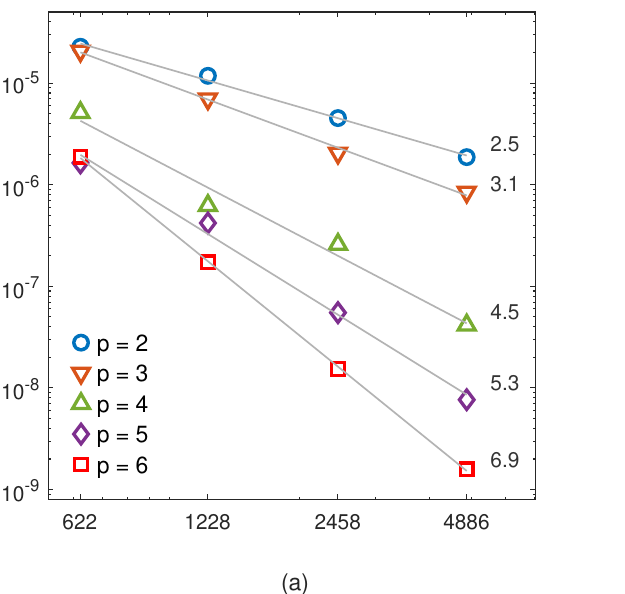}
\includegraphics[width=0.24\textwidth,clip=true,viewport=165 265 465 550]{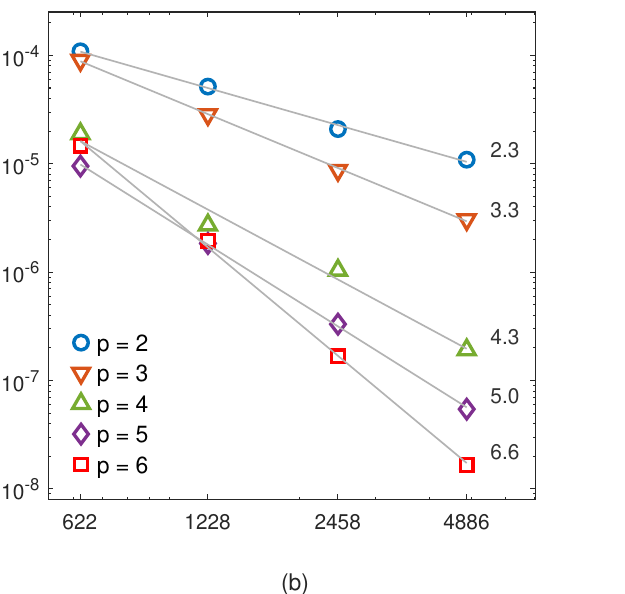}
\includegraphics[width=0.24\textwidth,clip=true,viewport=165 265 465 550]{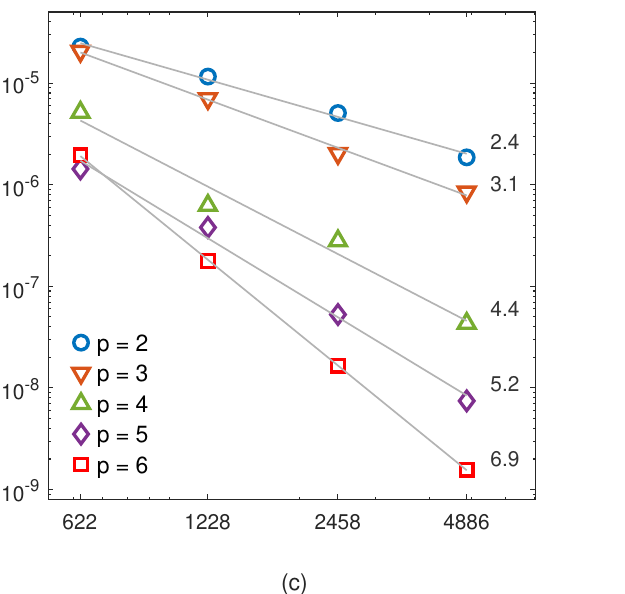}
\includegraphics[width=0.24\textwidth,clip=true,viewport=165 265 465 550]{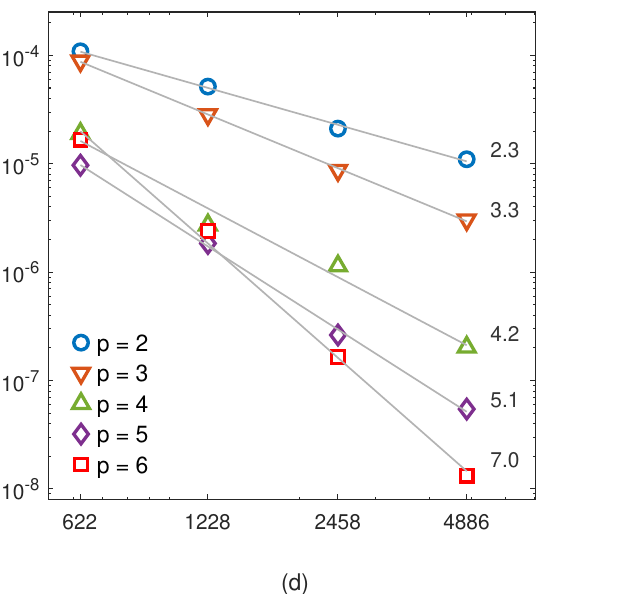}
\caption{\CapText Results for the test case 3.1 for $p=2,\ldots,6$. The $L^{1}$- and $L^{\infty}$- numerical solution errors for $u_1$ (panels (a) and (b), respectively) and $u_2$ (panels (c) and (d), respectively) are presented for a set of four grids with $622$, $1288$, $2458$, and $4886$ nodes. The corresponding
average convergence order is also presented for each $p$.}
\label{fig:Burgers_p100_convergencerror}
\end{figure}

We present the numerical solution error (using $L^{\infty}$- and $L^{1}$-norms) and the corresponding average convergence order (ACO) for $p=2,\ldots,6$ in Figure \ref{fig:Burgers_p100_convergencerror}.

Looking at ACO, we conclude that the optimal convergence order is obtained (or even exceeded) for both the $L^{\infty}$- and the $L^1$-errors, independently of the approximation degree. We also report that the ratio $L^{\infty}$-error / $L^1$-error obtained for $p=2,\ldots,6$ ($4.9 \pm 0.5$, $4.1 \pm 0.2$, $4.2 \pm 0.4$, $5.9 \pm 0.9$, and $9.8 \pm 2.0$, respectively) is well bounded, indicating the good quality of the numerical solution across the computational domain for all values of $p$. 

%------------
% Test case 2
%------------
\subsection{Test case 3.2}
\begin{figure}[ht]
\centering
\includegraphics[width=0.24\textwidth,clip=true,viewport=150 280 485 558]{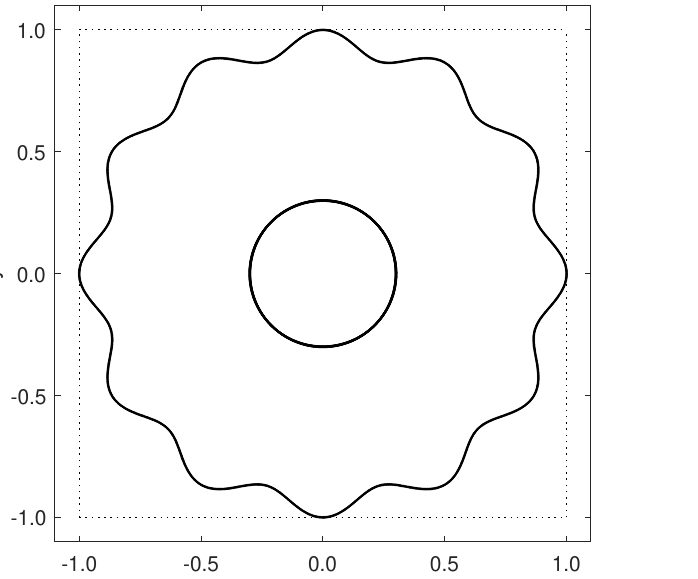}
\caption{\CapText Test case 3.2. Geometrical setting.}
\label{fig:test110geometry}
\end{figure}

We consider once again the two-dimensional Burgers' equations supplemented by Dirichlet boundary conditions over a non-polygonal domain expressed in polar coordinates with respect to the origin of the coordinate system (see Figure  \ref{fig:test110geometry}):
$$
\Omega =\left\{ (r,\theta)\in \mathbb{R}^{+} \times \left[ 0,2\pi\right[
: 0.3 \leq r \leq 0.95+0.05\sin(12\theta+\pi/2)\right\}.
$$

In order to perform a convergence analysis, we set $\bff$ and $\bu_D$
such that the unique solution of the BVP is given by
$$
u_1(r,\theta)=-q(r,\varepsilon)\sin\theta, \qquad u_2(r,\theta)=q(r,\varepsilon)\cos\theta,
$$
with
$$
q(r,\varepsilon )=\frac{r_4}{\varepsilon }\exp \left[ \frac{r_{3}\varepsilon }{%
\left( r-r_{1}\right) \left( r-r_{2}\right) }\right],  
$$
where we adopt $\varepsilon=10^{-2}$, while $r_1$, $r_2$, $r_3$, and $r_4$ are free parameters we shall set in the sequel. The solution of the BVP is invariant to rotation with respect to the origin of the coordinate system and we obtain
\begin{align*}
f_{1}( r,\theta) =-\frac{q^{2}}{r}\cos
\theta +\varepsilon \left( \frac{\partial ^{2}q}{\partial r^{2}}+\frac{1}{%
r\mathstrut }\frac{\partial q}{\partial r}-\frac{q}{r^{2}}\right) \sin
\theta,  \\
f_{2}(r,\theta) =-\frac{q^{2}}{r}\sin
\theta -\varepsilon \left( \frac{\partial ^{2}q}{\partial r^{2}}+\frac{1}{%
r\mathstrut }\frac{\partial q}{\partial r}-\frac{q}{r^{2}}\right) \cos
\theta.
\end{align*}

\begin{figure}[ht]
\centering
\includegraphics[width=0.24\textwidth,clip=true,viewport=135 280 470 558]{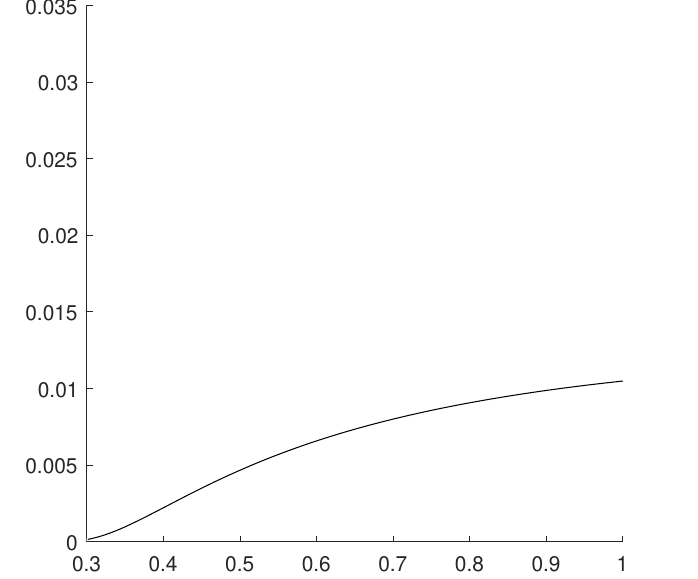}
\includegraphics[width=0.24\textwidth,clip=true,viewport=135 280 470 558]{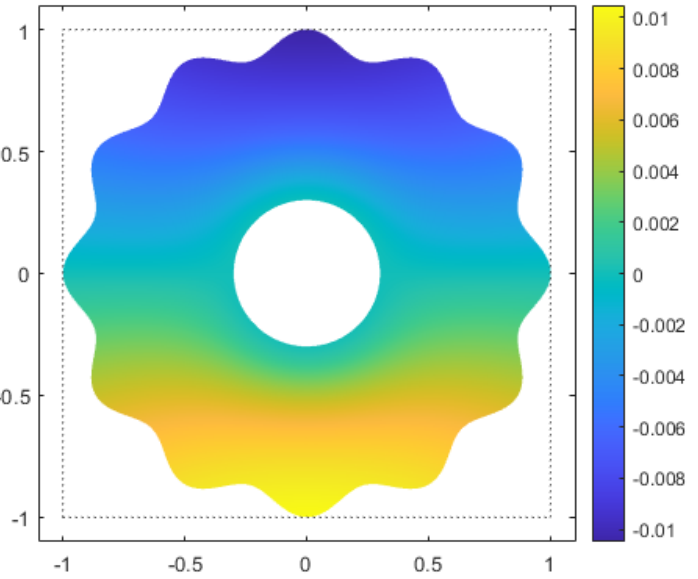}
\includegraphics[width=0.24\textwidth,clip=true,viewport=135 280 470 558]{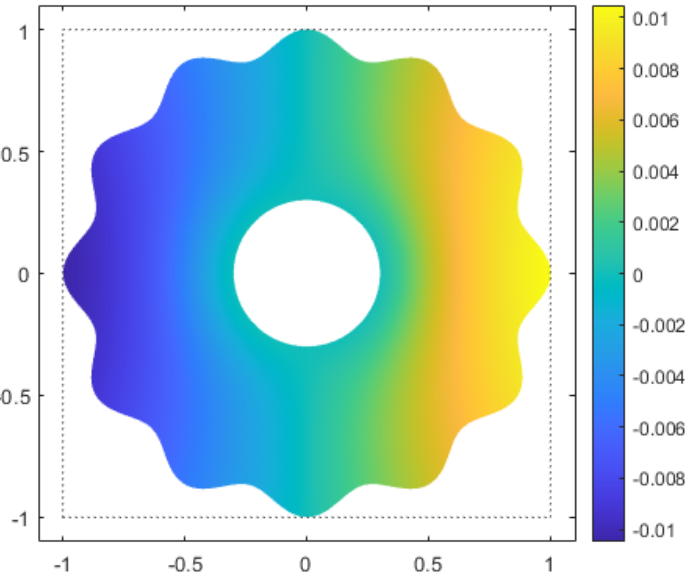}
\includegraphics[width=0.24\textwidth,clip=true,viewport=150 280 485 558]{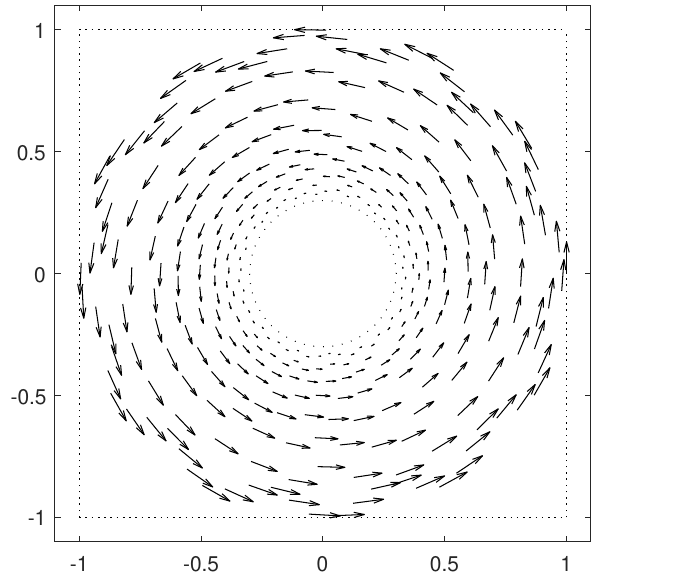}
\caption{\CapText Test case 3.2a. Shape of function $q$ (left panel) and exact solution (from left to right: $u_1$, $u_2$, $\bu)$.}
\label{fig:test110solution}
\end{figure}

\paragraph{Test case 3.2a}
We start the numerical simulations by setting $r_1=3$, $r_2=0.225$, $r_3=100$, and $r_4=2\times 10^{-4}$. The corresponding shape for the function $q$ (leading to a shallow slope curve) and the exact solution of the BVP are depicted in Figure \ref{fig:test110solution}.

\begin{figure}[ht]
\centering
\includegraphics[width=0.24\textwidth,clip=true,viewport=155 280 465 580]{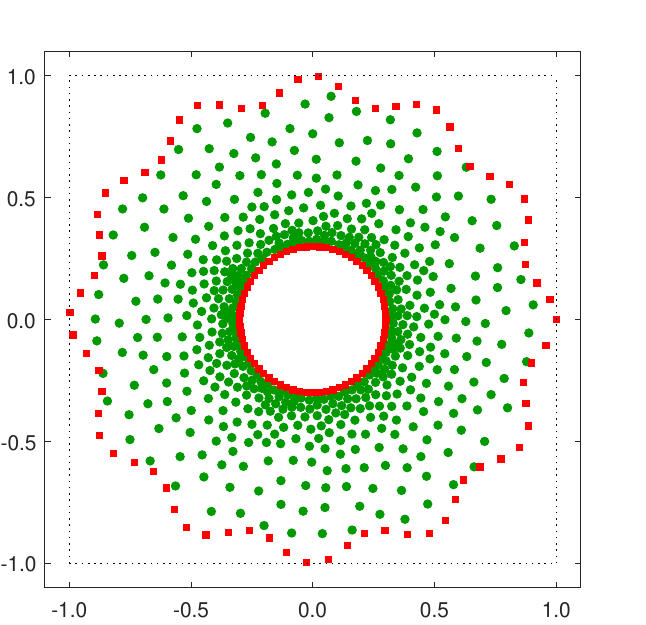}
\includegraphics[width=0.24\textwidth,clip=true,viewport=155 280 465 580]{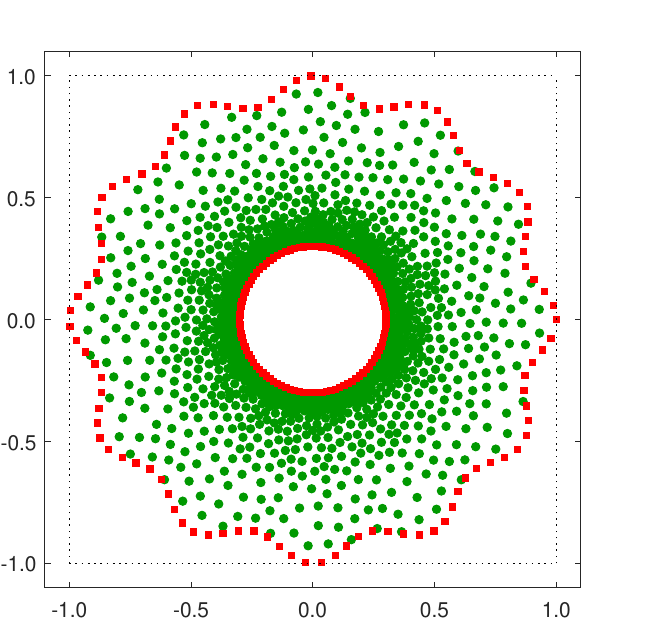}
\includegraphics[width=0.24\textwidth,clip=true,viewport=155 280 465 580]{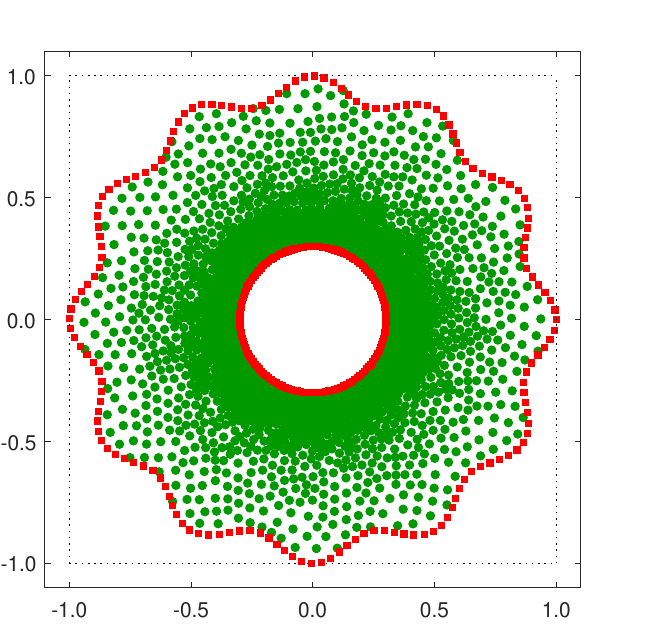}
\includegraphics[width=0.24\textwidth,clip=true,viewport=155 280 465 580]{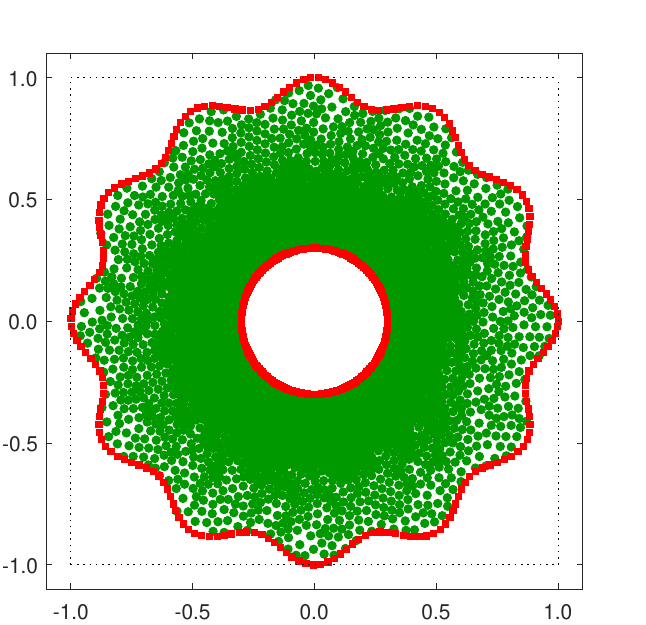}
\caption{\CapText Test case 3.2a. For $\overline{\Omega}$ the panels illustrate, from left to right, four grids obtained with a Delaunay triangulation procedure with $674$, $1316$, $2533$, and $5118$ nodes. The green filled circles and the red squares correspond, respectively, to interior nodes and boundary nodes where Dirichlet conditions are imposed.}
\label{fig:testcase110grid}
\end{figure}

For the numerical simulations, we use four grids with $674$, $1316$, $2533$, and $5118$ nodes, obtained with a Delaunay triangulation procedure. The number of nodes in the inner boundary doubles that of its outer counterpart. Figure \ref{fig:testcase110grid} illustrates the different grids used in the test case 3.2a numerical simulations. The nonlinearity is addressed using a fixed-point approach, as described previously.

%The results obtained for the error and convergence order of the numerical solution for $p=2,\ldots,6$ are presented in Tables \ref{tab:Burgers_p110_alpha=2} - \ref{tab:Burgers_p110_alpha=6}, respectively.

\begin{figure}[ht]
\centering
\includegraphics[width=0.24\textwidth,clip=true,viewport=165 265 465 550]{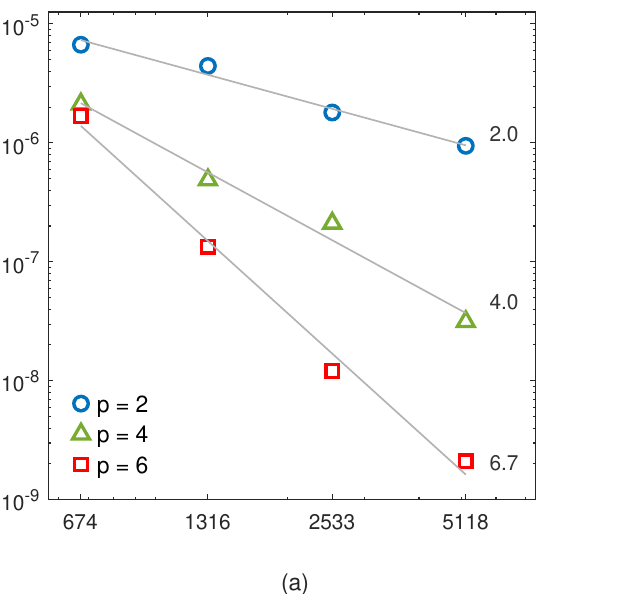}
\includegraphics[width=0.24\textwidth,clip=true,viewport=165 265 465 550]{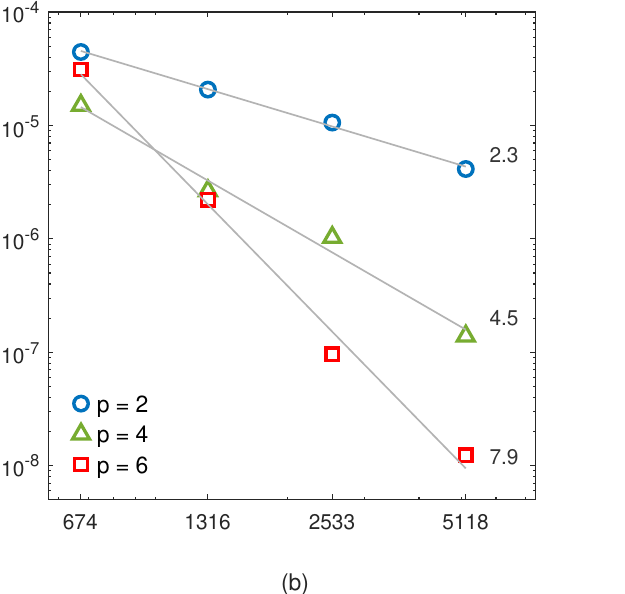}
\includegraphics[width=0.24\textwidth,clip=true,viewport=165 265 465 550]{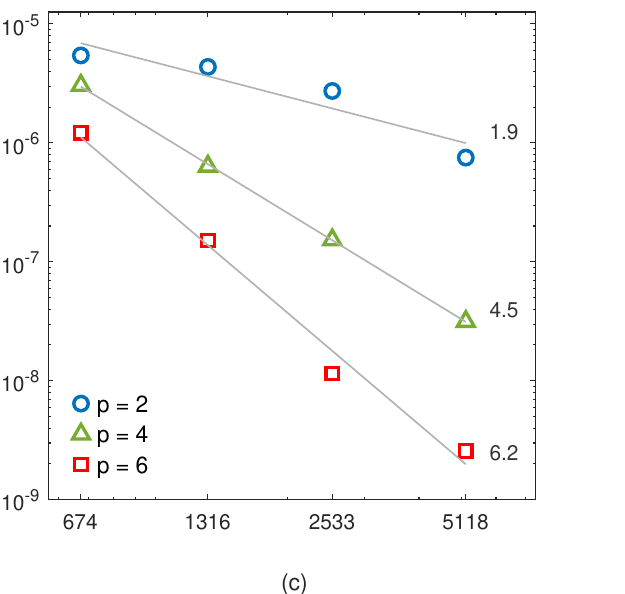}
\includegraphics[width=0.24\textwidth,clip=true,viewport=165 265 465 550]{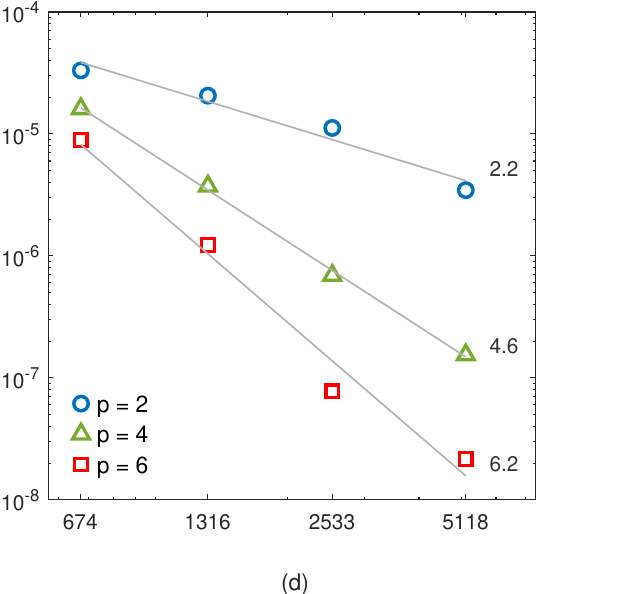}
\caption{\CapText Results for the test case 3.2a for $p=2,\ldots,6$ (the results obtained for $p=3$ and $p=5$ are not presented for the sake of clarity since they are very similar to those obtained for $p=2$ and $p=4$, respectively). The $L^{1}$- and $L^{\infty}$- numerical solution errors for $u_1$ (panels (a) and (b), respectively) and $u_2$ (panels (c) and (d), respectively) are presented for a set of four grids with $674$, $1316$, $2533$, and $5118$ nodes. The corresponding
average convergence order is also presented for each $p$.}
\label{fig:Burgers_p110_convergencerror}
\end{figure}

In Figure \ref{fig:Burgers_p110_convergencerror}, we present the numerical solution error (using $L^{\infty}$- and $L^{1}$-norms) and the corresponding ACO for $p=2,\ldots,6$.

From the analysis of the results obtained for the ACO, we conclude that the optimal convergence order is obtained (or even exceeded) for both the $L^{\infty}$- and the $L^1$-errors, as long as the approximation degree is even. When $p$ is odd, the ACO obtained equals the one obtained for the $p-1$ case, a feature already reported in our previous work.

We also report that the ratio $L^{\infty}$-error / $L^1$-error obtained for $p=2,\ldots,6$ ($5.1 \pm 0.8$, $3.2 \pm 0.2$, $5.3 \pm 0.7$, $4.0 \pm 0.7$, and $9.9 \pm 4.0$, respectively) remains moderate and very similar to the values obtained in case 3.1, indicating the good quality of the numerical solution for the whole computational domain for all values of $p$. 

\paragraph{Test case 3.2b}

\begin{figure}[ht]
\centering
\includegraphics[width=0.24\textwidth,clip=true,viewport=135 280 470 558]{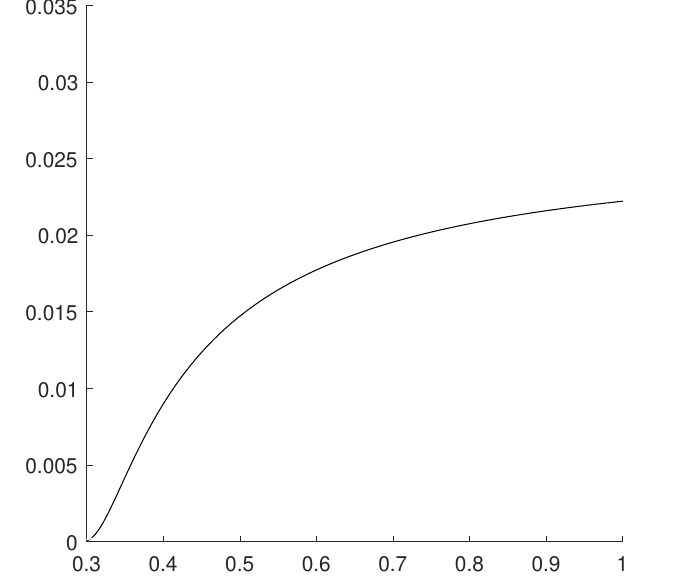}
\includegraphics[width=0.24\textwidth,clip=true,viewport=135 280 470 558]{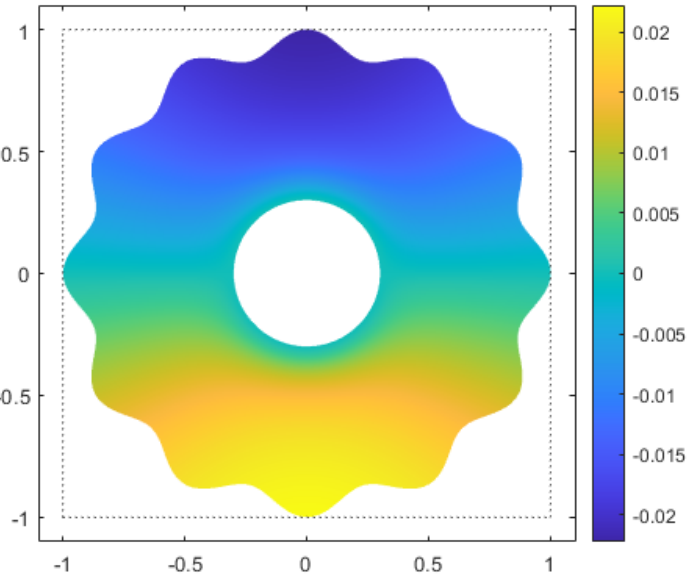}
\includegraphics[width=0.24\textwidth,clip=true,viewport=135 280 470 558]{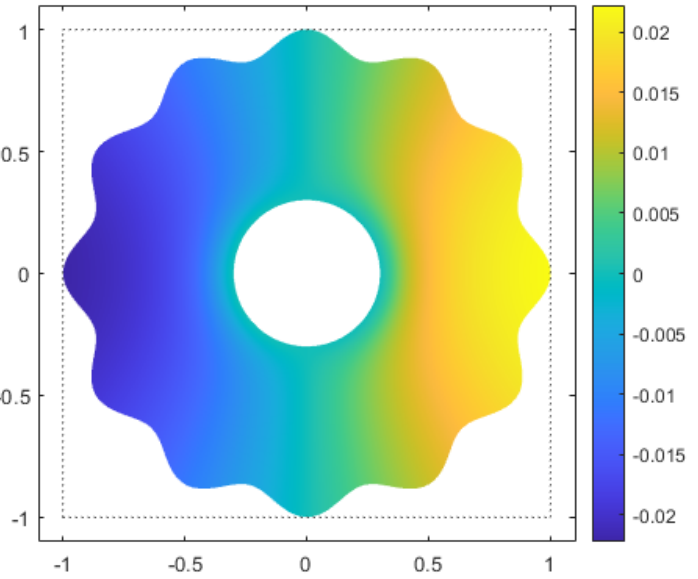}
\includegraphics[width=0.24\textwidth,clip=true,viewport=150 280 485 558]{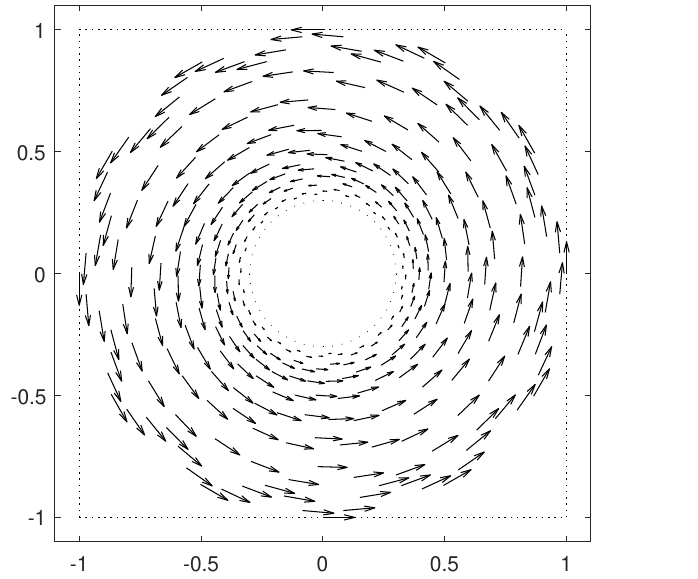}
\caption{\CapText Test case 3.2b. Shape of function $q$ (left panel) and exact solution (from left to right: $u_1$, $u_2$, $\bu$).}
\label{fig:test115solution}
\end{figure}

We now consider a solution that has an increased gradient near the inner boundary by setting $r_1=4$, $r_2=0.2775$, $r_3=50$, and $r_4=2.8\times 10^{-4}$. The corresponding shape of function $q$ and the exact solution of the BVP are illustrated in Figure \ref{fig:test115solution}.

\begin{figure}[ht]
\centering
\includegraphics[width=0.24\textwidth,clip=true,viewport=155 280 465 580]{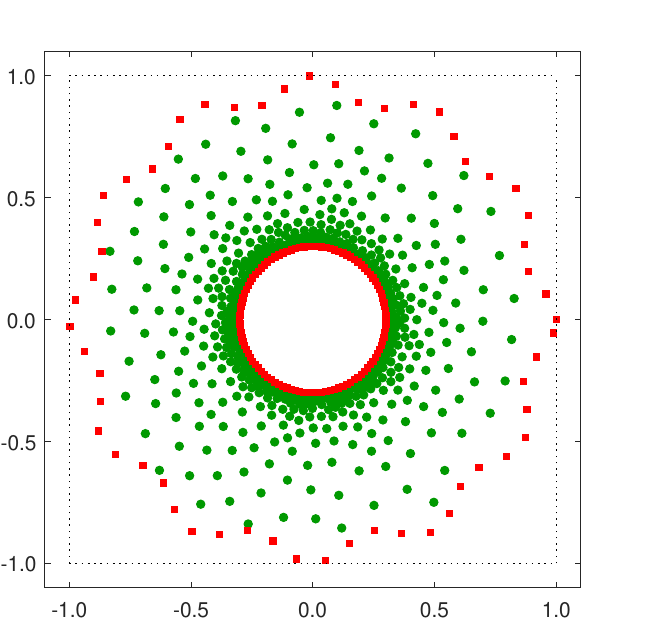}
\includegraphics[width=0.24\textwidth,clip=true,viewport=155 280 465 580]{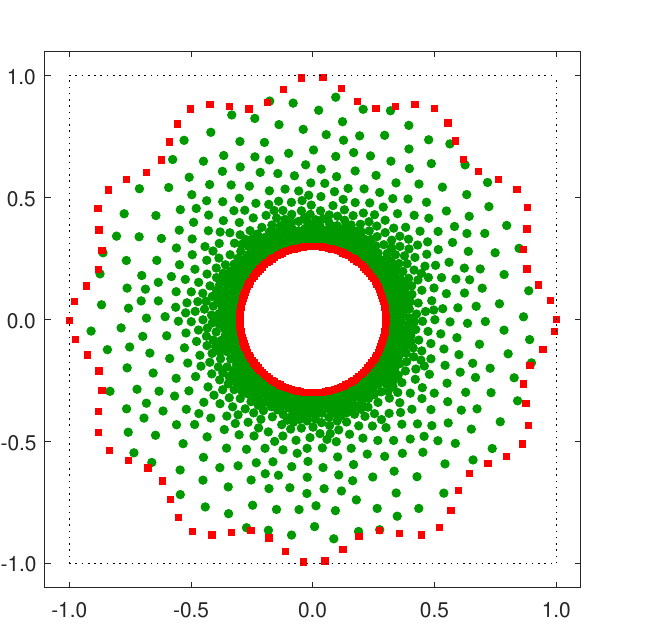}
\includegraphics[width=0.24\textwidth,clip=true,viewport=155 280 465 580]{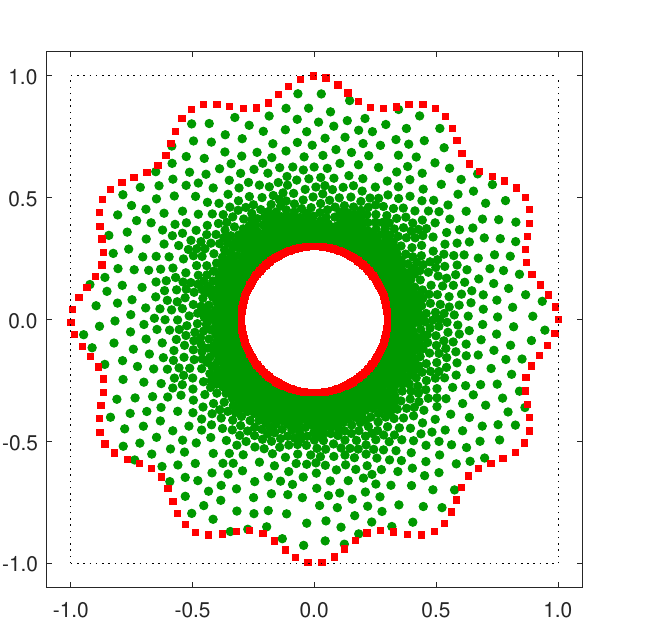}
\includegraphics[width=0.24\textwidth,clip=true,viewport=155 280 465 580]{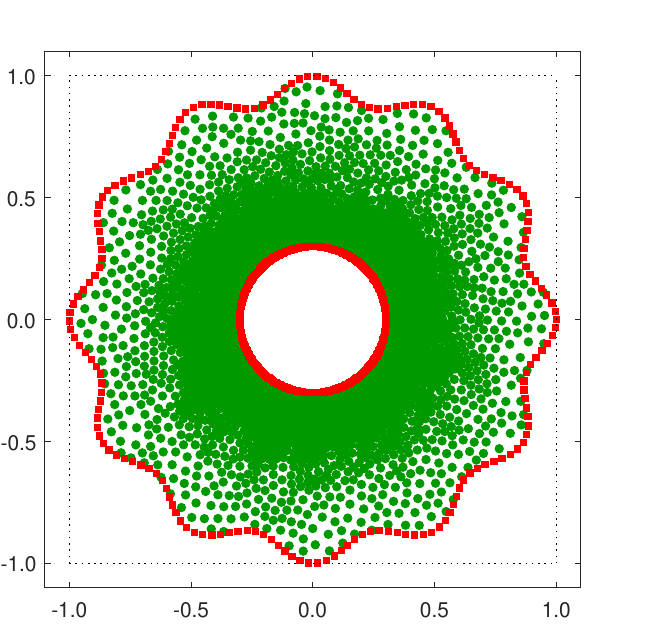}
\caption{\CapText Test case 3.2b. For $\overline{\Omega}$ the panels illustrate, from left to right, four grids obtained with a Delaunay triangulation procedure with $667$, $1317$, $2521$, and $5178$ nodes. The green filled circles and the red squares correspond, respectively, to interior nodes and boundary nodes where Dirichlet conditions are imposed.}
\label{fig:testcase115grid}
\end{figure}

For the numerical simulations, we use four grids with $667$, $1317$, $2521$, and $5178$ nodes, obtained with a Delaunay triangulation procedure. In order to account for the increased gradient near the inner boundary, the number of nodes in the inner boundary is now four times the number of nodes in its outer counterpart. We note that the total number of nodes in the four grids is very similar to the one in test case 3.2a. Figure \ref{fig:testcase115grid} illustrates the different grids used in the test case 3.2b simulations.

\begin{figure}[ht]
\centering
\includegraphics[width=0.24\textwidth,clip=true,viewport=165 265 465 550]{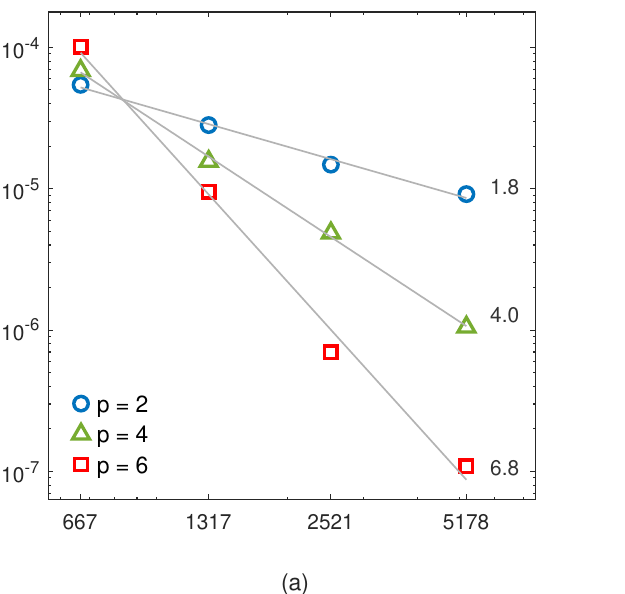}
\includegraphics[width=0.24\textwidth,clip=true,viewport=165 265 465 550]{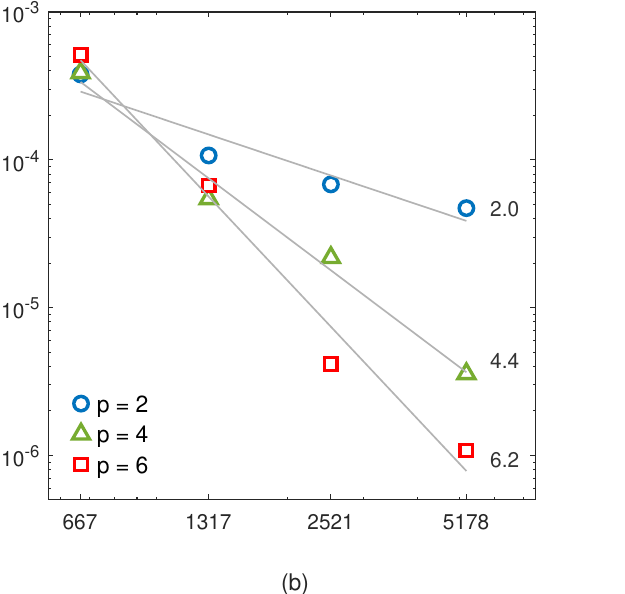}
\includegraphics[width=0.24\textwidth,clip=true,viewport=165 265 465 550]{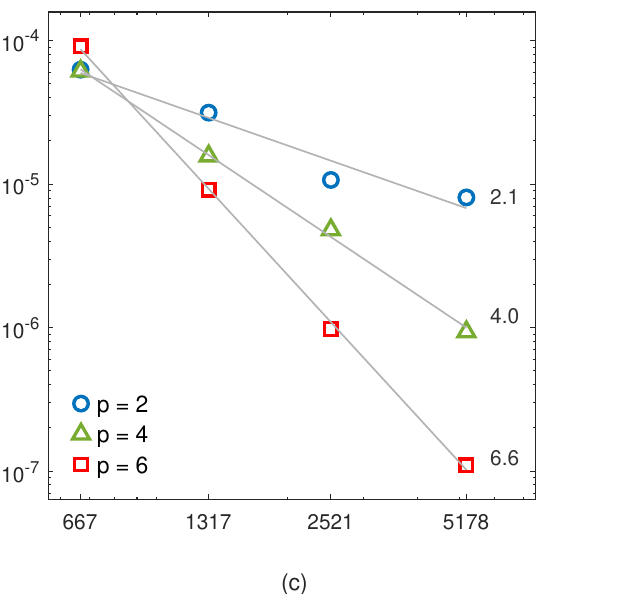}
\includegraphics[width=0.24\textwidth,clip=true,viewport=165 265 465 550]{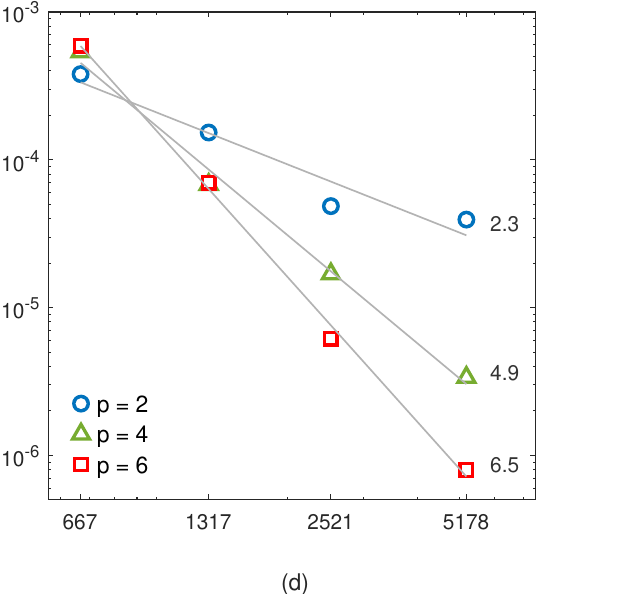}
\caption{\CapText Results for the test case 3.2b for $p=2,\ldots,6$ (the results obtained for $p=3$ and $p=5$ are not presented for the sake of clarity since they are very similar to those obtained for $p=2$ and $p=4$, respectively). The $L^{1}$- and $L^{\infty}$- numerical solution errors for $u_1$ (panels (a) and (b), respectively) and $u_2$ (panels (c) and (d), respectively) are presented for a set of four grids with $667$, $1317$, $2521$, and $5178$ nodes. The corresponding
average convergence order is also presented for each $p$.}
\label{fig:Burgers_p115_convergencerror}
\end{figure}

The results obtained for the numerical solution error and convergence order for $p=2,\ldots,6$ are presented in Figure \ref{fig:Burgers_p115_convergencerror}.

We conclude that even though the errors are higher than those obtained in test case 3.2a, increasing the gradient of function $q$ near the inner boundary does not impact the convergence order, since it follows the same trends already described for test case 3.2a. Likewise, the results observed for the ratio $L^{\infty}$-error / $L^1$-error for $p=2,\ldots,6$ ($5.1 \pm 0.8$, $4.3 \pm 0.5$, $4.7 \pm 1.5$, $4.9 \pm 1.1$, $7.0 \pm 1.2$, respectively) are also in line with the good values obtained for the previous test cases.

\paragraph{Test case 3.2c}

\begin{figure}[ht]
\centering
\includegraphics[width=0.24\textwidth,clip=true,viewport=135 280 470 558]{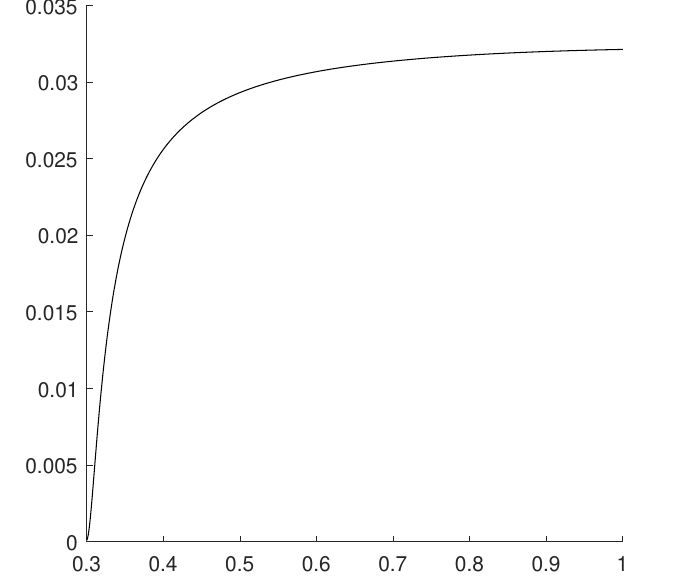}
\includegraphics[width=0.24\textwidth,clip=true,viewport=135 280 470 558]{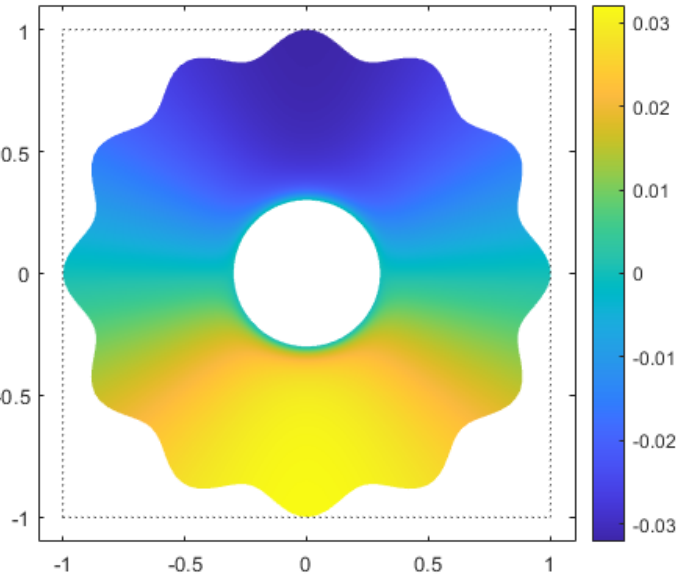}
\includegraphics[width=0.24\textwidth,clip=true,viewport=135 280 470 558]{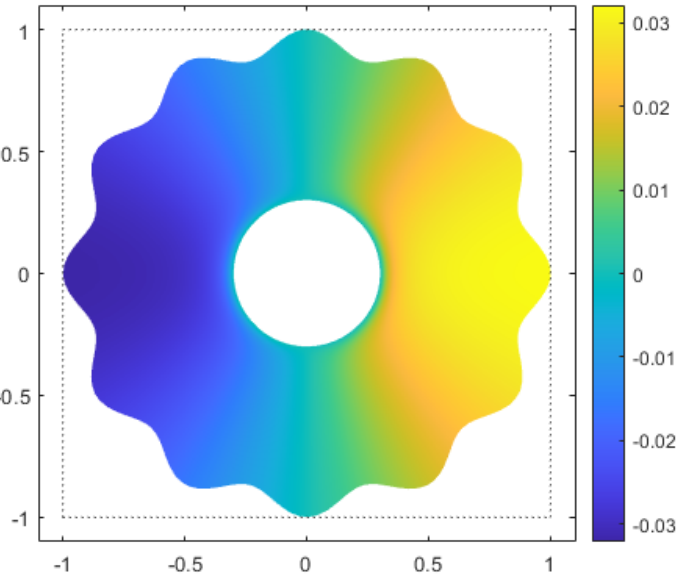}
\includegraphics[width=0.24\textwidth,clip=true,viewport=150 280 485 558]{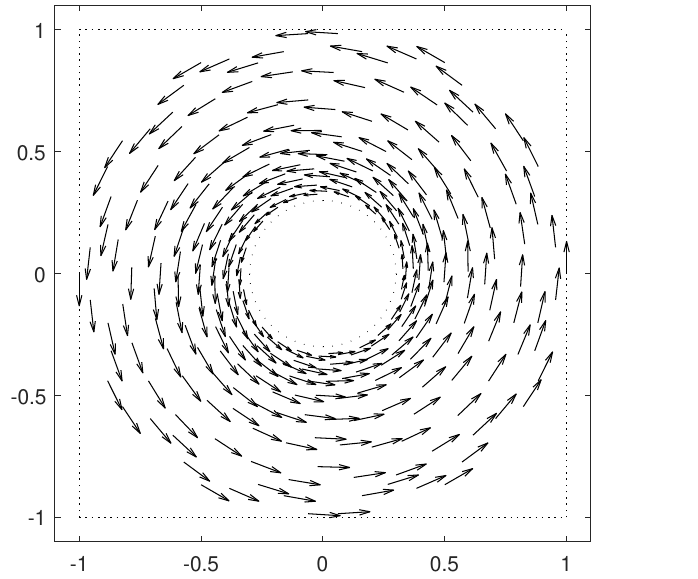}
\caption{\CapText Test case 3.2c. Shape of function $q$ (left panel) and exact solution (from left to right: $u_1$, $u_2$, $\bu$).}
\label{fig:test120solution}
\end{figure}

Finally, setting $r_1=2$, $r_2=0.295$, $r_3=5$, and $r_4=3.45\times 10^{-4}$, we consider a solution that is even steeper near the inner boundary, mimicking a boundary layer scenario. The corresponding shape for function $q$ and the exact solution of the BVP are shown in Figure \ref{fig:test120solution}.

\begin{figure}[ht]
\centering
\includegraphics[width=0.24\textwidth,clip=true,viewport=155 280 465 580]{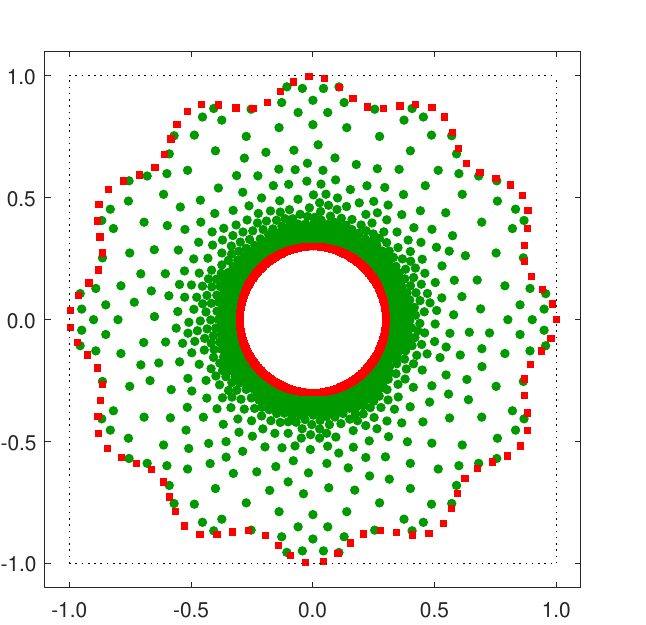}
\includegraphics[width=0.24\textwidth,clip=true,viewport=155 280 465 580]{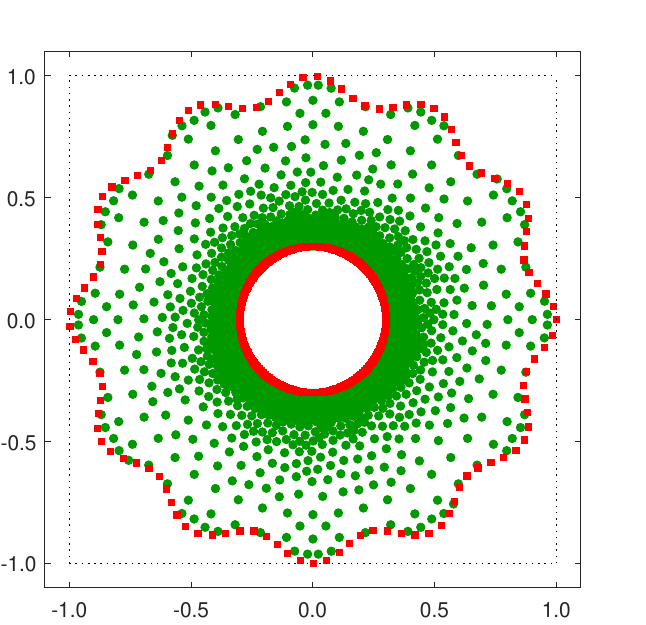}
\includegraphics[width=0.24\textwidth,clip=true,viewport=155 280 465 580]{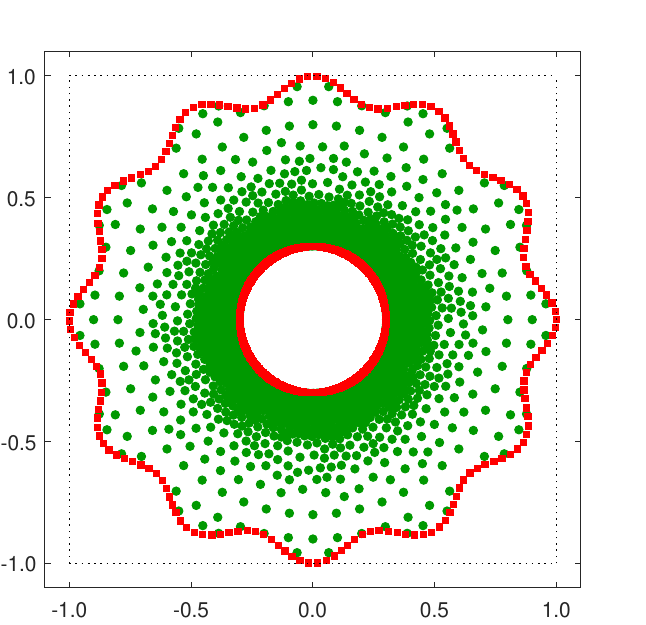}
\includegraphics[width=0.24\textwidth,clip=true,viewport=155 280 465 580]{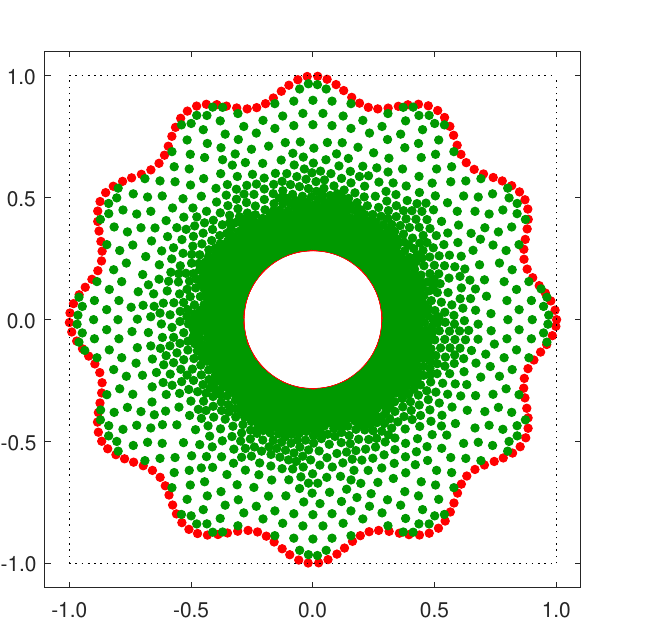}
\caption{\CapText Test case 3.2c. For $\overline{\Omega}$ the panels illustrate, from left to right, four grids obtained with a Delaunay triangulation procedure with $2711$, $3723$, $5159$ and $7169$ nodes. The green filled circles and the red squares correspond, respectively, to interior nodes and boundary nodes where Dirichlet conditions are imposed.}
\label{fig:testcase120grid}
\end{figure}

In order to account for the rather steep gradient near the inner boundary, the number of nodes in the inner boundary is set to ten times the number of nodes in its outer counterpart. The density of nodes in the outer region has to be kept similar to the one used in test case 3.2b in order to ensure the convergence of the fixed point method for all values of $p$. Therefore, the number of nodes has to be increased with respect to the test cases 3.2a and 3.2b. Hence, for the numerical simulations, we use four grids with $2711$, $3723$, $5159$, and $7169$ nodes, obtained with a Delaunay triangulation procedure. Figure \ref{fig:testcase120grid} illustrates the four grids used in the test case 3.2c simulations. 

\begin{figure}[ht]
\centering
\includegraphics[width=0.24\textwidth,clip=true,viewport=165 265 465 550]{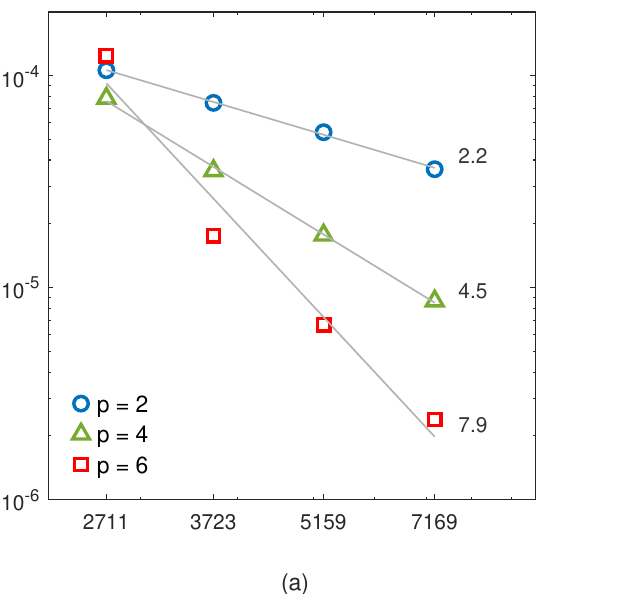}
\includegraphics[width=0.24\textwidth,clip=true,viewport=165 265 465 550]{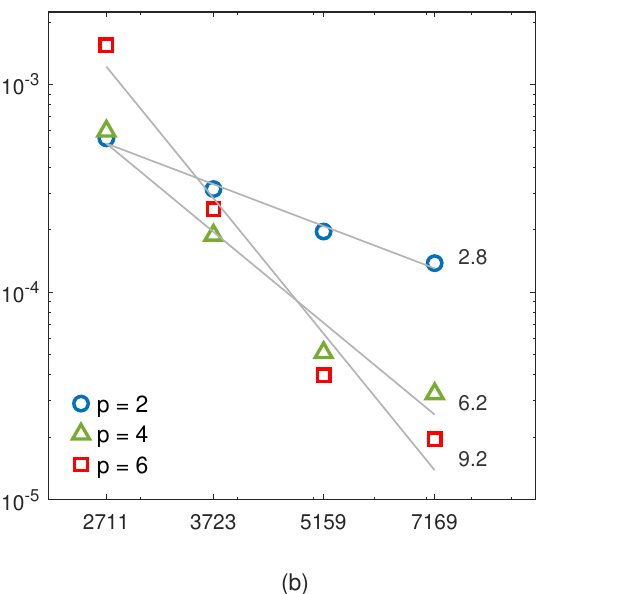}
\includegraphics[width=0.24\textwidth,clip=true,viewport=165 265 465 550]{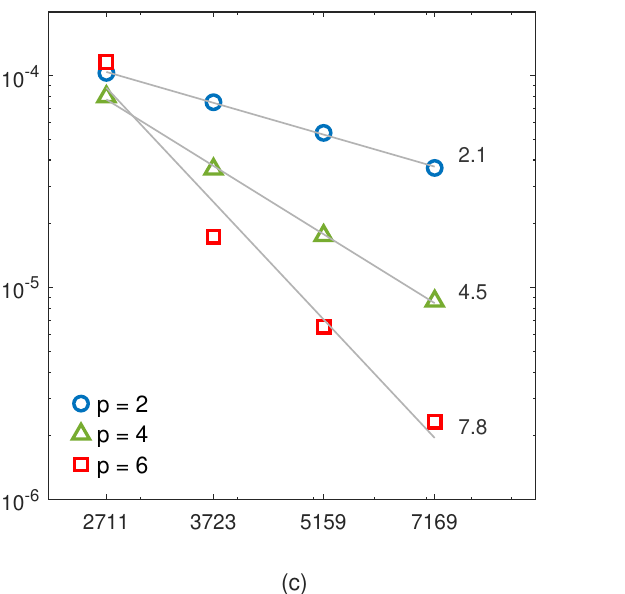}
\includegraphics[width=0.24\textwidth,clip=true,viewport=165 265 465 550]{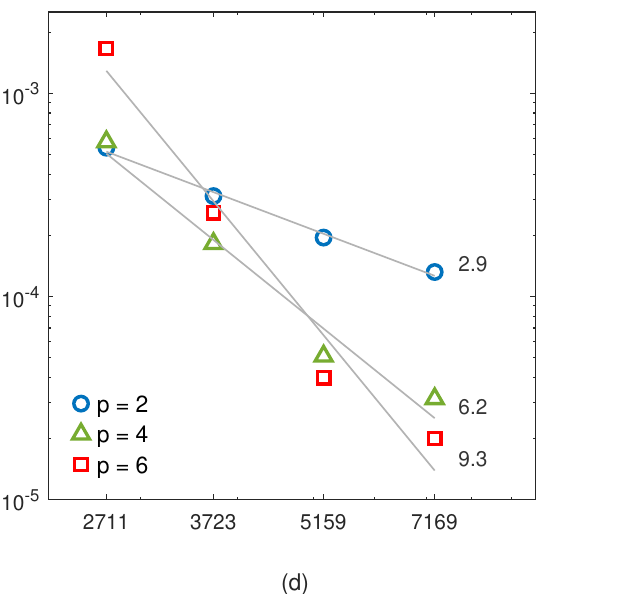}
\caption{\CapText Results for the test case 3.2c for $p=2,\ldots,6$ (the results obtained for $p=3$ and $p=5$ are not presented for the sake of clarity since they are very similar to those obtained for $p=2$ and $p=4$, respectively). The $L^{1}$- and $L^{\infty}$- numerical solution errors for $u_1$ (panels (a) and (b), respectively) and $u_2$ (panels (c) and (d), respectively) are presented for a set of four grids with $2711$, $3723$, $5159$, and $7169$ nodes. The corresponding
average convergence order is also presented for each $p$.}
\label{fig:Burgers_p120_convergencerror}
\end{figure}

The results obtained for the numerical solution error and the ACO for $p=2,\ldots,6$ are presented in Figure \ref{fig:Burgers_p120_convergencerror}.

The results show that the use of a more demanding function $q$ inevitably leads to higher errors, but the convergence order remains essentially the same as in test cases 3.2a and 3.2b, reflecting the good ability of the numerical scheme to handle a boundary layer-like scenario. On the other hand, the results observed for the ratio $L^{\infty}$-error / $L^1$-error for $p=2,\ldots,6$ ($4.2 \pm 0.6$, $4.1 \pm 0.4$, $4.8 \pm 1.6$, $3.2 \pm 0.5$, $10.6 \pm 3.2$, respectively) support the good quality of the numerical solution across the computational domain for all values of $p$.

%--------------
% Test case 3.3
%--------------

\subsection{Test case 3.3}

For this last test case concerning the two-dimensional Burgers' equations, we aim to test the ability of the numerical scheme to address Neumann boundary conditions in the presence of complex geometries. As in the previous test case, we take $\varepsilon=10^{-2}$. We now consider a 2D nozzle-like domain given by

\begin{equation*}
\Omega =\left\{ \bx\in \mathbb{R}^{2}:0.5<x_1<1,\ -h(x_1)<x_2<h(x_1)\right\} ,
\end{equation*}
where
\begin{align*}
\Gamma_{t}& =\left\{ \bx\in \mathbb{R}^{2}:0.5<x_1<1,\ x_2=h(x_1)\right\} , \\
\Gamma_{b}& =\left\{ \bx\in \mathbb{R}^{2}:0.5<x_1<1,\ x_2=-h(x_1)\right\} , \\
\Gamma_{l}& =\left\{ \bx\in \mathbb{R}^{2}:x_1=0.5,\ -h(x_1)\leq x_2\leq h(x_1)\right\} , \\
\Gamma_{r}& =\left\{ \bx\in \mathbb{R}^{2}:x_1=1,\ -h(x_1)\leq x_2\leq h(x_1)\right\} ,
\end{align*}

with $h(x_1)=0.025\arccos(0.1\cos40x_1)$. The solution of the BVP satisfies homogeneous Neumann conditions on the upper and 
lower boundaries, $\Gamma_{t}$ and $\Gamma_{b}$, and Dirichlet conditions on the left and right boundaries, $\Gamma_{l}$
and $\Gamma_{r}$, that is,
\begin{equation*}
\begin{array}{l}
\bu(\bx)=\bu_{D}(\bx),\text{\quad on } \Gamma_{l}\cup \Gamma _{r},\smallskip \\ 
\partial _{\mathbf{n}}{\bu}(\bx) = \mathbf{0},\text{\quad on } \Gamma _{t}\cup \Gamma _{b},
\end{array}
\end{equation*}
where $\mathbf{n}$ is the outer unit normal vector. Therefore, we set $\bmu(\bx_i)=\mathbf{1}$ and $\bnu(\bx_i)=\mathbf{0}$, for all $\bx_i \in\Gamma_{l} \cup \Gamma_{r}$, and $\bmu(\bx_i)=\mathbf{0}$ and $\bnu(\bx_i)=\mathbf{1}$, for all $\bx_i \in\Gamma_{t} \cup \Gamma_{b}$. In order to perform a convergence analysis we set $\bff$ and $\bu_D$ such that
\begin{equation*}
u_1(\bx)=\frac{1}{5}(1+\cos40x_1\cos40x_2),\text{\quad }u_2(\bx)=\frac{1}{5}(1-\cos40x_1\cos40x_2),
\end{equation*}
is the unique solution of the BVP. The choice of function $h$ ensures homogeneous Neumann conditions hold on $\Gamma _{t}$ and 
$\Gamma_{b}$ for this particular choice of the solution of the BVP (\textit{cf.} \cite{CLP21}).

\begin{figure}[ht]
\centering
\includegraphics[width=0.2275\textwidth,clip=true,viewport=140 280 480 558]{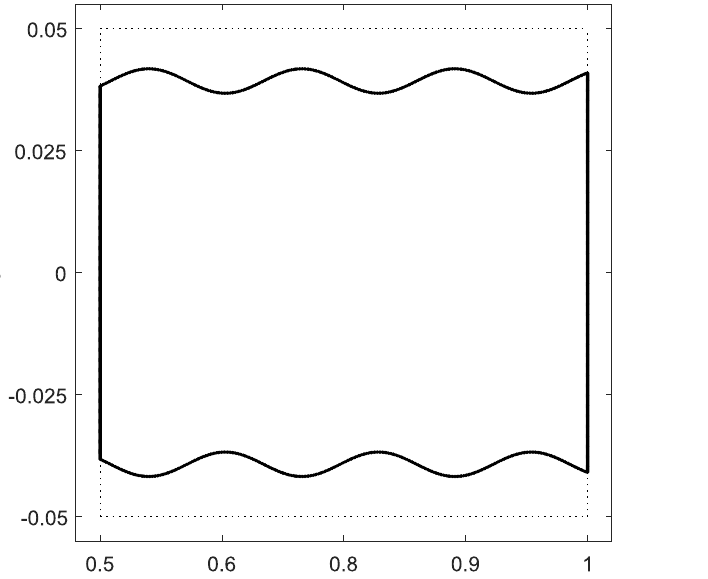}
\includegraphics[width=0.24\textwidth,clip=true,viewport=120 280 480 558]{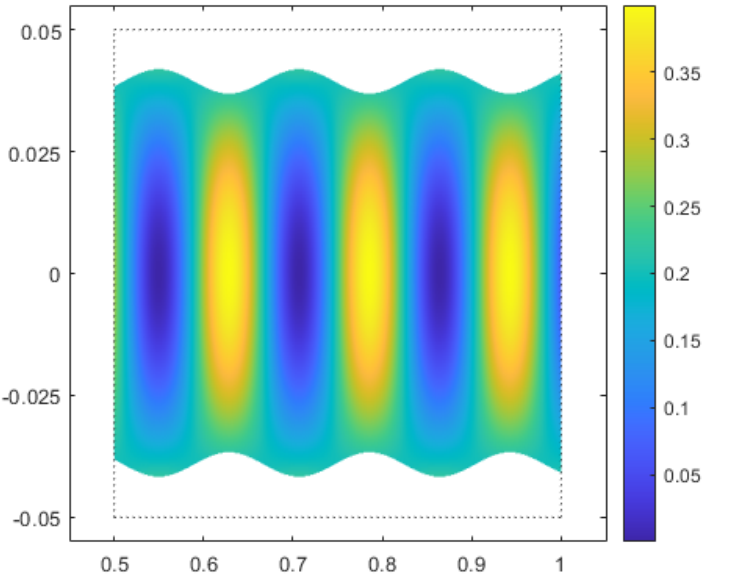}
\includegraphics[width=0.24\textwidth,clip=true,viewport=120 280 480 558]{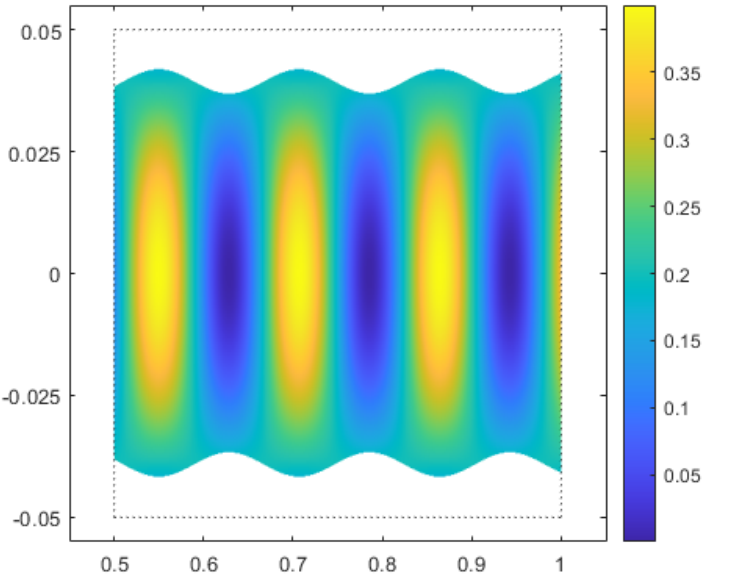}
\includegraphics[width=0.24\textwidth,clip=true,viewport=120 280 480 558]{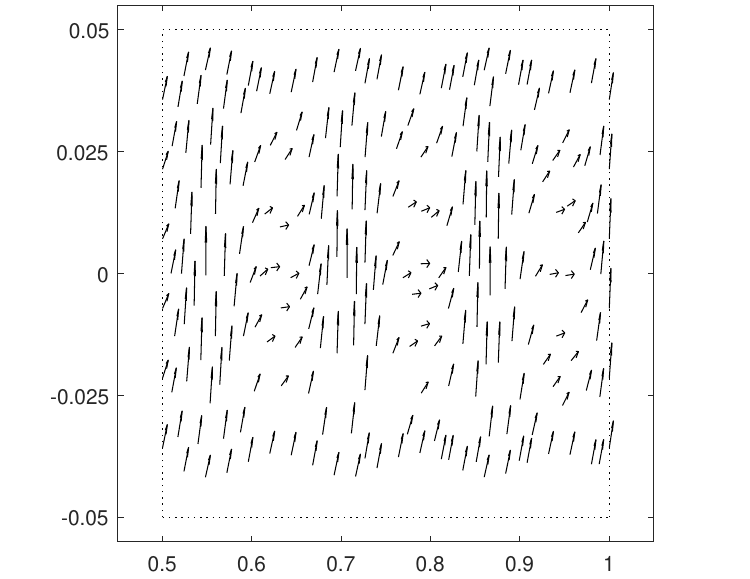}
\caption{\CapText Test case 3.3. Geometrical setting (left panel) and exact solution (from left to right: $u_1$, $u_2$, $\bu$).}
\label{fig:test130solution}
\end{figure}

The geometrical setting and the exact solution of the BVP are illustrated in Figure \ref{fig:test130solution}.

\begin{figure}[ht]
\centering
\includegraphics[width=0.24\textwidth,clip=true,viewport=145 280 465 580]{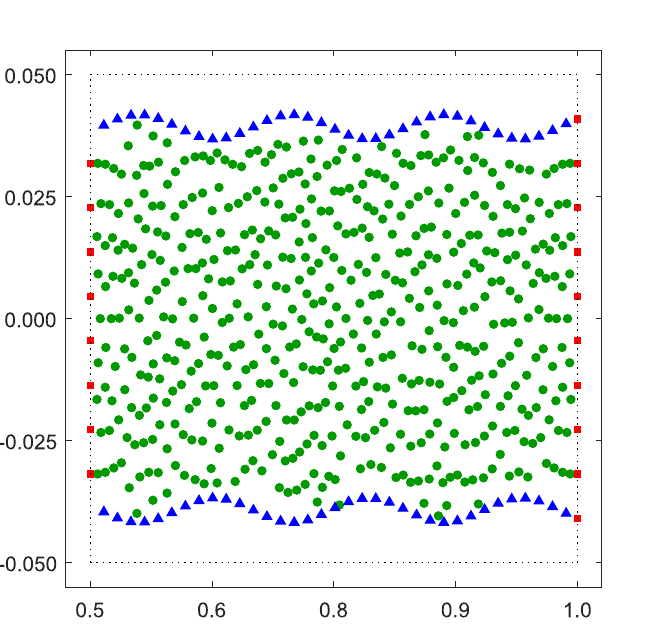}
\includegraphics[width=0.24\textwidth,clip=true,viewport=145 280 465 580]{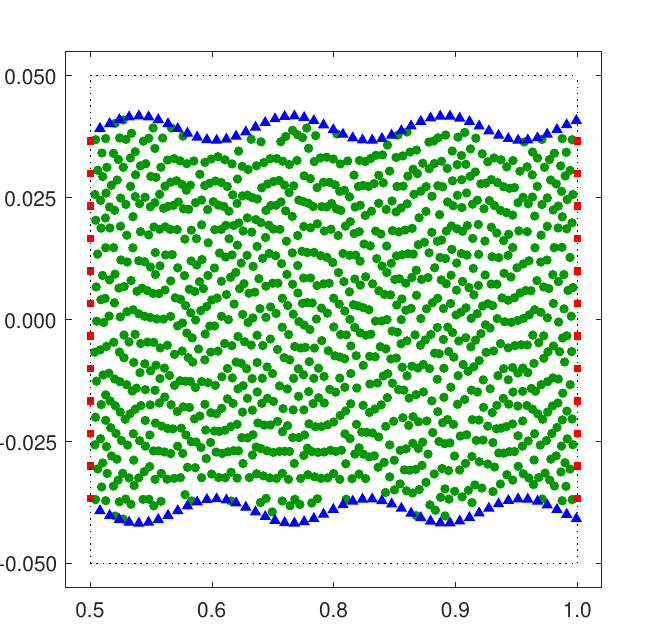}
\includegraphics[width=0.24\textwidth,clip=true,viewport=145 280 465 580]{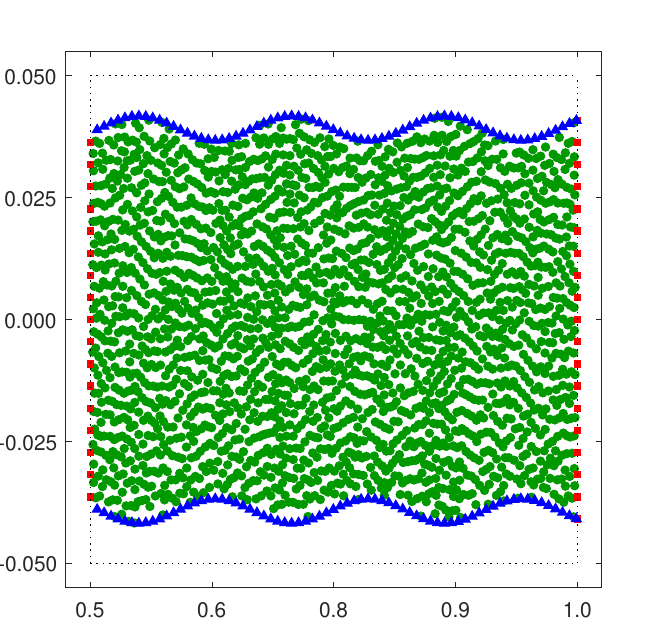}
\includegraphics[width=0.24\textwidth,clip=true,viewport=145 280 465 580]{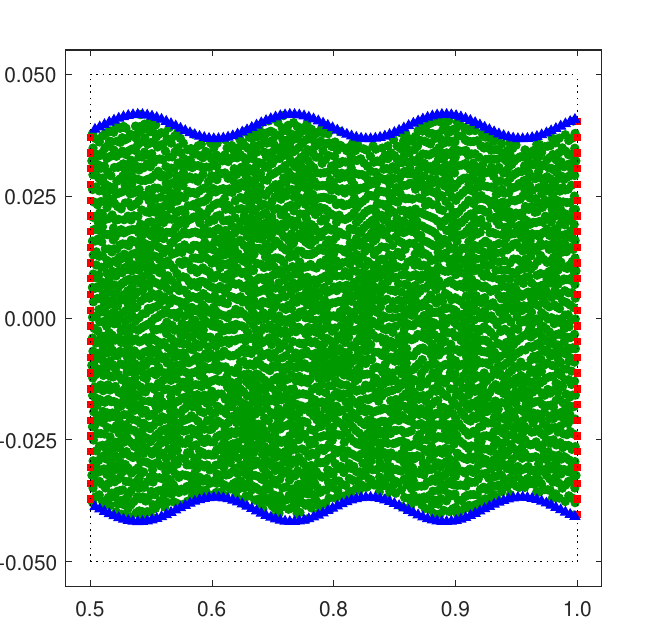}
\caption{\CapText Test case 3.3. For $\overline{\Omega}$ the panels illustrate, from left to right, four grids obtained with a Delaunay triangulation procedure with $625$, $1264$, $2505$, and $5008$ nodes. The green filled circles, the red squares and the blue triangles correspond, respectively, to interior nodes, boundary nodes where Dirichlet conditions are imposed and boundary nodes where Neumann conditions hold.}
\label{fig:testcase130grid}
\end{figure}

For the numerical simulations, we use four grids with $625$, $1264$, $2505$, and $5008$ nodes, obtained using a Delaunay triangulation procedure. Figure \ref{fig:testcase130grid} illustrates the different grids used in test case 3.3 numerical simulations.

\begin{figure}[ht]
\centering
\includegraphics[width=0.24\textwidth,clip=true,viewport=165 265 465 550]{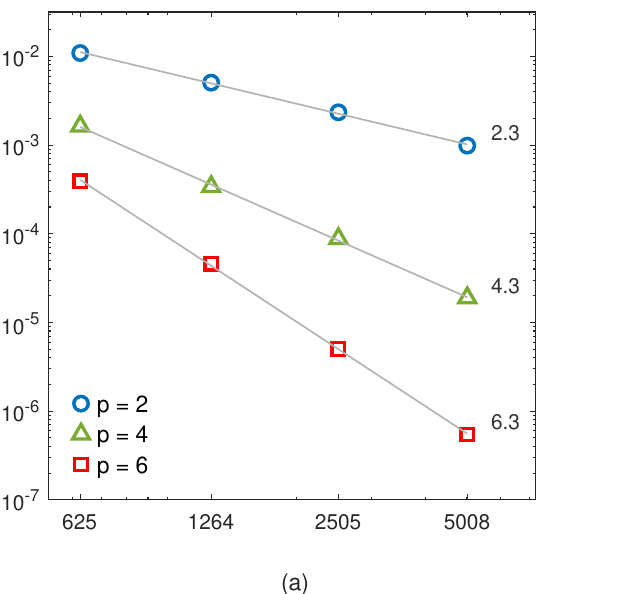}
\includegraphics[width=0.24\textwidth,clip=true,viewport=165 265 465 550]{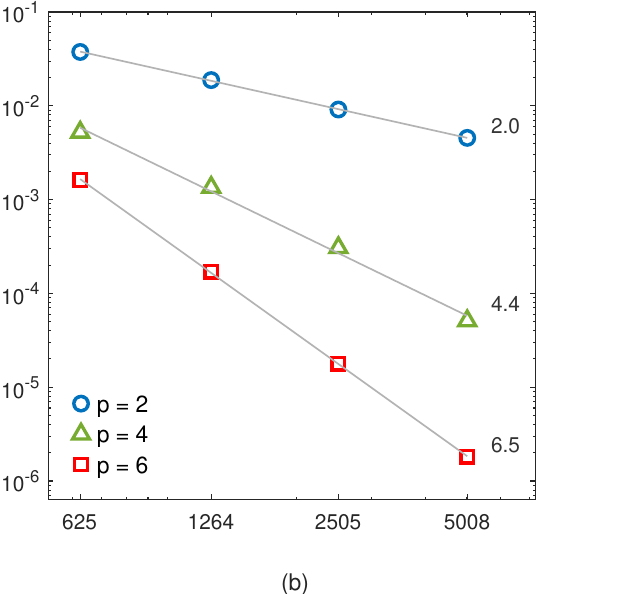}
\caption{\CapText Results for the test case 3.3 for $p=2,\ldots,6$ (the results obtained for $p=3$ and $p=5$ are not presented for the sake of clarity since they are very similar to those obtained for $p=2$ and $p=4$, respectively). The $L^{1}$- and $L^{\infty}$- numerical solution errors for $u_1$ (panels (a) and (b), respectively) are presented for a set of four grids with $625$, $1264$, $2505$, and $5008$ nodes. The corresponding
average convergence order is also presented for each $p$. We note that the results obtained for $u_2$ are omitted since they are identical to those presented for $u_!$.}
\label{fig:Burgers_p130_convergencerror}
\end{figure}

The results obtained for the numerical solution error and convergence order for $p=2,\ldots,6$ are presented in Figure \ref{fig:Burgers_p130_convergencerror}.

From the analysis of the results, we conclude that optimal convergence orders are obtained for the two norms used to evaluate the numerical solution error when the polynomial approximation degree is even. As in test case 3.2, the results obtained for $p=3$ and $p=5$, in particular the ACO, are rather similar to those found for $p=2$ and $p=4$, respectively. The good quality of the numerical solution across the computational domain is also confirmed by the low values obtained for the ratio $L^{\infty}$-error / $L^1$-error ($3.9 \pm 0.8$, $2.6 \pm 0.3$, $3.4 \pm 0.9$, $3.1 \pm 1.1$, $3.7 \pm 0.5$, for $p=2,\ldots,6$, respectively), which remains below $4$.

\section{Conclusion}
The general finite difference method, by its nature, enables one to wisely choose the nodes of the stencil attached to a generic node $i$ while other popular meshless methods, such as FPsM or FPtM, always consider the nodes belonging to the compact ball given by a user-dependent radius. Such a degree of freedom provides both accuracy and stability while maintaining a small stencil size. As a continuation of our work on the scalar linear PDE, we prove that the strategy works very well for the non-linear steady-state problem using a fixed-point technique based on a local linearization. Only low local Péclet/Reynolds-number discrete problems were considered to assess the optimization algorithm initially based on an isotropic stencil. The next stage will be to introduce some upwind bias, both for the stencil and the weights, to tackle convection-dominated problems and even hyperbolic systems with smooth solutions.

\section*{Acknowledgments}
\noindent The authors Stéphane Clain and Jorge Figueiredo were funded by national funds from FCT - Fundação para a Ciência e a Tecnologia, I.P., within the context of the project HiPMos4C setting with reference 2023.16854.ICDT2023.16854.ICDT, https://doi.org/10.54499/2023.16854.ICDT.\\
\noindent Stéphane Clain would like to acknowledge the support of Centro de Matemática of the Coimbra University (CMUC), funded by the Portuguese Foundation for Science and Technology, under the references  UIDB/00324/2020 (https://doi.org/10.54499/UIDB/00324/2020) and UIDP/00324/2020 (https://doi.org/10.54499/UIDP/00324/2020).\\
\noindent Jorge Figueiredo was supported by national funds through FCT  - Fundação para a Ciência e a Tecnologia, under the Program Contract of the Centre of Mathematics of the University of Minho (CMAT/UM), UID/00013/2025 (https://doi.org/10.54499/UID/00013/2025).

\section*{Conflicts of Interest} The authors declare that there are no conflicts of interest.

\bibliographystyle{elsarticle-num}
\bibliography{references.bib}

\end{document}